\documentclass[11pt]{article}
\usepackage{amsmath, amsfonts, amssymb, mathtools}
\usepackage{graphicx, amsthm,color}
\usepackage{cite}
\usepackage[english]{babel}
\usepackage[utf8]{inputenc}
\usepackage[T1]{fontenc}
\usepackage{amscd}
\usepackage{newlfont}
\usepackage{enumerate}
\usepackage{float}
\usepackage{makecell} 
\usepackage{graphicx}
\usepackage{multirow,booktabs}
\usepackage{caption, subcaption}
\usepackage{cancel}
\usepackage{tikz-cd}

\newcommand{\R}{\mathbb{R}}
\newcommand{\N}{\mathbb{N}}
\newcommand{\e}{\varepsilon}
\theoremstyle{definition}
\usepackage{thmtools}
\declaretheoremstyle[%
  headfont=\bfseries,%
  headpunct={.},%
  bodyfont=\normalfont,
  notefont=\bfseries
]{definitionstyle}
\declaretheoremstyle[%
  headfont=\bfseries,%
  headpunct={.},%
  bodyfont=\normalfont\itshape,
  notefont=\bfseries
]{theoremstyle}
\declaretheoremstyle[%
  headfont=\itshape,%
  headpunct={.},%
  bodyfont=\normalfont,
  notefont=\bfseries
]{remarkstyle}
\declaretheorem[style=definitionstyle,name=Definition,numberwithin=section]{definition}
\declaretheorem[style=theoremstyle,name=Theorem,sibling=definition]{theorem}

\declaretheorem[style=theoremstyle,name=Proposition,sibling=definition]{proposition}
\declaretheorem[style=theoremstyle,name=Lemma,sibling=definition]{lemma}
\declaretheorem[style=definitionstyle,name=Example,sibling=definition]{example}
\declaretheorem[style=remarkstyle,name=Remark,sibling=definition]{remark}
\usepackage{caption, subcaption}
\declaretheoremstyle[
  spaceabove=1mm,
  spacebelow=1mm,
  headfont=\bfseries,
  headpunct={.},
  notefont=\bfseries
]{algorithmstyle}
\declaretheorem[style=algorithmstyle,name=Assumption,sibling=definition]{assumption}
\usepackage{hyperref}
\usepackage{cleveref}
\crefname{theorem}{theorem}{theorems}
\crefname{definition}{definition}{definitions}
\crefname{corollary}{corollary}{corollaries}
\crefname{proposition}{proposition}{propositions}
\crefname{remark}{remark}{remarks}
\crefname{example}{example}{examples}
\crefname{lemma}{lemma}{lemmas}
\crefname{assumption}{assumption}{assumptions}
\crefname{subsection}{subsection}{subsections}
\numberwithin{equation}{section}

\newcounter{saveenum}

\begin{document}

\begin{center}
{\LARGE \textbf{Solving Nonlinear Absolute Value Equations}}
\end{center}
\smallskip

\begin{center}
{\Large \textsc{Aris Daniilidis, Mounir Haddou, Tr\'i Minh L\^e  \bigskip \\ Olivier Ley, Phi Hoang Tran}}
\end{center}

\bigskip

\noindent\textbf{Abstract.} In this work, we show that several problems naturally represented as Nonlinear Absolute Value Equations (NAVE) can be reformulated as Nonlinear Complementarity Problems (NCP) and efficiently solved using smoothing regularization techniques under mild assumptions. As far as we know, this is the first numerical approach that directly deals with NAVE. We also identify a technical assumption commonly utilized in smoothing techniques and prove its equivalence to a classical {\L}ojasiewicz inequality at infinity, validating its non-restrictive nature. Furthermore, we extend established error estimates for NCP solvers to derive error bounds for NAVE problems under weaker assumptions. We illustrate the effectiveness of our approach through applications including asymmetric ridge optimization and nonlinear ordinary differential equations.

\bigskip


\vspace{0.6cm}

\textit{Secondary}: 15B48, 65K15, 90C59.

\section{Introduction}\label{sec_introduction}
In the last two decades, absolute value equation problems (in short, \ref{AVE} problems) have been extensively studied in the literature. This interest is justified by the fact that this class of problems already covers a wide spectrum of applications: indeed, numerous problems stemming from real-life applications, as for instance all mixed integer linear programming problems, can be reformulated as \ref{AVE} problems. It is also well-known, see \cite{P2009,M2014} \textit{e.g.}, that \ref{AVE} problems admit an equivalent description as Linear Complementarity problems (in short, \ref{LCP}). The exact definitions of an \ref{AVE} problem and a \ref{LCP} are recalled below. Dealing efficiently with these problems is thus paramount. 

\medskip

The literature on this subject contains several theoretical results for existence as well as conditions guaranteeing uniqueness of the solution \cite{M2006, RHF2014,R2009,L2013}.  
Concurrently, there are also various numerical approaches to solve an \ref{AVE} problem. Generally speaking, these methods can be divided into at least three categories  \cite{ACT2023}: iterative linear algebra based methods (also known as projective methods), semi-smooth Newton-like methods and smoothing methods. 
The aforementioned methods generally require some assumption on the matrix involved in the \ref{AVE} problem. 
In particular, the classes of $P_0$-matrices and $P$-matrices (recalled below) turn out to be relevant in this study \cite{AHM2018}. 

\medskip
 
In this work, we consider a natural generalization of (linear) \ref{AVE} problems to nonlinear ones, known as Nonlinear Absolute Value Equations (in short, \ref{AVE}). This more general framework encompasses new applications including ridge regression models, bounded constrained nonlinear systems of equations, and  stiff Ordinary Differential Equations (in short, stiff ODE). This approach to deal with the aforementioned problems, based on \ref{NAVE}, is to the best of our knowledge, completely new in the literature.

\medskip
Our main contribution is as follows: first, we demonstrate that, similar to the connection between an \ref{AVE} problem and an \ref{LCP}, a \ref{NAVE} can be associated with a nonlinear complementarity problem \ref{NCP}. Although this associated \ref{NCP} is generally implicit, we show that it can be handled efficiently. Leveraging the extensive literature on the existence, uniqueness, and numerical resolution of \ref{NCP} (e.g., \cite{HH2010, M2006, RHF2014, W2020}), we propose a new method for solving the \ref{NAVE} problem, based on the smoothing technique developed in \cite{HM2014} and further extended in \cite{OHBA2022}. The proposed approach is explained in \Cref{subsec_tranformation}, while in \Cref{sec_applications} we discuss applications.  \smallskip\newline
To ease  the reading we start with some definitions and settings. The \textit{Absolute Value Equation
problem} \ref{AVE} is defined as follows:
\begin{equation}
\text{find }x\in\mathbb{R}^{d}:\qquad Ax-|x|=b, \tag{AVE}\label{AVE}
\end{equation}
where $A$ is a $(d\times d)$-matrix and $b\in\mathbb{R}^{d}$.
Throughout this work, given
$x=(x^{1},\cdots,x^{d})^{T}\in\mathbb{R}^{d}$, 
we use the notation $|x|:=(|x^{1}|,\cdots,|x^{d}|)^{T}$ componentwise to denote a vector in $\mathbb{R}^d_+$.\smallskip\newline
Denoting by $I$ the identity matrix of $\mathbb{R}^{d}$ and assuming that either
$A-I$ or $A+I$ is invertible, the \ref{AVE} problem can be transformed to a \textit{Linear
Complementarity Problem} \ref{LCP}. Indeed, setting (coordinate by coordinate)
\begin{equation}
\left\{
\begin{array}
[c]{c}
y=x^{+}=\max\,\{\phantom{-}x,0\}\smallskip\\
z=x^{-}=\max\,\{-x,0\}
\end{array}
\right.  \label{eq:pos-neg}
\end{equation}
and performing the transformation $x=y-z$ and $|x|=y+z$ we obtain a \ref{LCP} problem
\[
\left(  A-I\right)  y-\left(  A+I\right)  z=b.
\]
Therefore for
\begin{equation}
\left\{
\begin{array}
[c]{l}
M:=(A+I)^{-1}(A-I)\smallskip\\
q:=(A+I)^{-1}(-b)
\end{array}
\right.  \qquad\text{or respectively\qquad}\left\{
\begin{array}
[c]{l}
\widetilde{M}:=(A-I)^{-1}(A+I)\smallskip\\
\widetilde{q}:=(A-I)^{-1}b
\end{array}
\right.  \label{mounir1}
\end{equation}
we obtain a \ref{LCP} problem
\begin{equation}
  \left\{
\begin{array}
[c]{l}
z=My+q\smallskip\\
0\leq y\perp z\geq0
\end{array}
\right.  \qquad\text{or respectively\qquad}\left\{
\begin{array}
[c]{l}
y=\widetilde{M}z+\widetilde{q}\smallskip\\
0\leq y\perp z\geq0
\end{array}.
\right.  \tag{LCP}\label{LCP}
\end{equation}
The notation \(0 \le y \perp z \ge 0\) denotes the standard complementarity
condition, that is, 
\[
y, z \in \R^d_+ \qquad \text{and} \qquad y^j z^j = 0, \quad \text{ for every $j \in \{1, \cdots, d \}$}.
\]
\ref{LCP} can be solved provided that $M$ (respectively $\widetilde{M}$) is a $P$-matrix (see below for details). It is important to
notice that this property can be traced back to the matrix $A$; in particular, the property is ensured whenever the singular values of $A$ are all greater than~$1$. Notice that this condition guarantees invertibility of both $A-I$ and $A+I$. 

\medskip

Solving \ref{LCP} under the assumption that $M$ (respectively $\widetilde{M}$) is a $P$-matrix has been treated in several works (see \cite{BG}). In this case, it can be shown
that the problem has a unique solution $(\bar{y},\bar{z})$ yielding that
$\bar{x}:=\bar{y}-\bar{z}$ is the (unique) solution of \ref{AVE}. Moreover, this solution can be obtained numerically, via smoothing regularization
techniques (see~\cite{AHM2018, HM2014, OHBA2022} and \Cref{subsec_smoothing_techniques} below).\smallskip\newline
{Inspired by the fruitful link between \ref{AVE} and \ref{LCP},
we propose a new method of solving
the following \textit{Nonlinear Absolute Value Equation} \ref{NAVE}
\begin{equation}
\text{Find }x\in\mathbb{R}^{d}:\qquad F(x)-|x|=0\quad\text{(coordinatewise),}
\tag{NAVE}\label{NAVE}
\end{equation}
where $F:\mathbb{R}^{d}\rightarrow\mathbb{R}^{d}$ is a (nonlinear) mapping. By
introducing new variables $y=x^{+}$ and $z=x^{-}$ (\emph{cf.}~\eqref{eq:pos-neg}) so that $x=y-z$ and $|x|=y+z$, \ref{NAVE}
becomes
\begin{eqnarray}\label{eq-Fyz}
F(y-z)-(y+z)=0.
\end{eqnarray}
By setting $z=H(y)$ (which is possible under regularity assumptions on $F$, see Lemma~\ref{implicite-fct}), we can transform \ref{NAVE} into a \textit{Nonlinear
Complementarity Problem} \ref{NCP}:
\begin{equation}
\left\{
\begin{array}
[c]{l}
H(y)=F(y-H(y))-y\smallskip\\
0\leq y\perp H(y)\geq0.
\end{array}
\right.  \tag{NCP}\label{NCP}
\end{equation}
As we shall see  in \Cref{subsec_smoothing_techniques}, even though the function $H$ is only defined implicitly, it is still possible to solve~\ref{NCP} numerically provided we are able to guarantee that $H$ is a $P_{0}$-map (see Definition~\ref{def_Poper}), a notion which generalizes $P_0$-matrices (\textit{c.f.} Lemma~\ref{guarant-Pmap}).

\medskip
To our knowledge, no existing works directly address \ref{NAVE}, despite the formulation of various practical applications in this form. A significant contribution of this work is the characterization of a common technical assumption on smoothing functions for \ref{NCP}; we prove in Theorem~\ref{thm-mounir} that this assumption is equivalent to a classical {\L}ojasiewicz inequality at infinity, demonstrating that this assumption is natural and verifiable.

\medskip

This paper is organized as follows. In \Cref{sec_preliminary_and_problem}, we describe our approach for reformulating the \ref{NAVE} problem as a \ref{NCP} problem. 
We then examine the 
uniqueness of solutions to \ref{NAVE}, followed by an in-depth analysis of the implicit mappings arising from the reformulation. 
Next, we introduce smoothing techniques for solving \ref{NCP} problems and adapt them to address \ref{NAVE} problems in \Cref{sec_method}.
The second subsection of \Cref{sec_method} outlines the complete solver for \ref{NAVE} and then establishes its associated error estimates. 
In \Cref{sec_applications}, we demonstrate the effectiveness of the proposed algorithm on various applications, including asymmetric ridge regression, sparse optimization, and stiff ODE.
The paper concludes with a summary of contributions and potential directions for future research.


\section{Preliminaries and setting of the problem}
\label{sec_preliminary_and_problem}

\subsection{Definitions and preliminaries}

Given a $(d\times d)$ matrix $A$ and $I\subset \{1,2,\cdots,d\}$,  we denote by $A_{II}$ the submatrix made up of the rows and columns of $I$.

\begin{definition}[$P_0$-matrix and $P$-matrix]\label{def_Pmatrix} 
A matrix $A$ is called a $P_{0}$-matrix (respectively, $P$-matrix)
if one of the following equivalent properties holds
 \begin{enumerate}[(i).]
  \item for every $I\subset \{1,2,\cdots,d\}$, $\text{det}(A_{II}) \geq 0$ (respectively $\text{det}(A_{II}) > 0$);
  \item
for every $x=(x^1,\cdots , x^d)^T \in \R^d$, $x\not= 0$, 
     \begin{eqnarray*}
      \max_{\substack{1 \le j \le d \\ x^j \ne 0}} (Ax)^j x^j \geq 0 
       \qquad
       \Big(\text{respectively } 
       \max_{\substack{1 \le j \le d}} (Ax)^j x^j > 0 \Big);
      \end{eqnarray*}
   \item for every $I\subset \{1,2,\cdots,d\}$, the real eigenvalues of $A_{II}$ are nonnegative (resp. strictly positive).
 \end{enumerate}
\end{definition}
The proof for the equivalence of the above definitions can be found in \cite[Sections 3.3 and 3.4]{CPS2009}.
The notion of $P$-matrix can be generalized to general nonlinear maps
$H:\mathbb{R}^{d}\rightarrow\mathbb{R}^{d}$ as follows.

\begin{definition}\label{def_Poper}
A map $H = (H^1, \cdots, H^d)^T : \mathbb{R}^d \to \mathbb{R}^d$ is said to be 
a $P_0$-\textit{map} (resp. $P$-\textit{map}) if for every $x, y \in \mathbb{R}^d$ with $x \ne y$, it holds
    \begin{align*}
        \max\limits_{\substack{1 \leq i \leq d \\ x^i \ne y^i}} (x^i - y^i)(H(x)^i - H(y)^i) \geq 0 \quad \Bigg(\text{resp. } \max\limits_{\substack{1 \leq i \leq d}} (x^i - y^i)(H(x)^i - H(y)^i) > 0 \Bigg).
    \end{align*}
\end{definition}
\noindent

We refer the reader to~\cite{MR73} for further results about $P_0$- and $P$-maps.
To conclude, we state the following technical lemmas, which are needed to establish the transformation between \ref{NAVE} and \ref{NCP}.


\begin{lemma}[{\hspace*{-1mm}\cite[Corollary 5.3, Theorem 5.8]{MR73}}]  \label{lem-map-mat}
Let $H:\mathbb{R}^{d}\rightarrow \mathbb{R}^{d}$ be $C^1$. 
Then $H$ is a $P_0$-map if and only if, for every $x\in\R^d$, the Jacobian matrix $\nabla H(x)$ is a $P_0$-matrix.
\end{lemma}



\begin{lemma}[A characterization of $P_{0}$-matrices] \label{lem-P0-invert}
Let $A$ be a $(d\times d)$ matrix. Then $A$ is a $P_0$-matrix if and only if, for every diagonal matrix $\Delta_1$ with strictly positive entries and for every nonnegative diagonal matrix $\Delta_2$, the matrix $\Delta_1 +\Delta_2 A$ is invertible.
\end{lemma}
\noindent See \Cref{proof_lem-P0-invert} for the proof of this lemma.
\subsection[]{Transforming a \ref{NAVE} into a \ref{NCP}}\label{subsec_tranformation}

The following lemma gives conditions under which~\eqref{eq-Fyz} can be written as a \ref{NCP} as discussed at the end of \Cref{sec_introduction}, by setting $z=H(y)$ or $y=\widetilde{H}(z)$ for some suitable maps $H$ or $\widetilde{H}$.

\begin{lemma}\label{implicite-fct}
Assume that the mapping $F$ in \ref{NAVE} is $C^1$ in a neighborhood of the point ${x_*}=y_*-z_*\in\R^d$ which is assumed to be a solution of~\eqref{eq-Fyz}. Then it holds:
\begin{itemize}    
\item[$(i)$.] If $F-I$ is a $P_0$-map, then there exists a $C^1$ map $H:\mathbb{R}^d\rightarrow\mathbb{R}^d$ defined in a neighborhood of~$y_*$ such that $z_*=H(y_*)$
and $y_*$ is a solution to the following \ref{NCP} problem:
\begin{equation}
\left\{
\begin{array}
[c]{l}
H(y)=F(y-H(y))-y\smallskip\\
0\leq y\perp H(y)\geq0.
\end{array}
\right.  \label{eqNCP1}
\end{equation}

\item[$(ii)$.] If $-(F+I)$ is a $P_0$-map, then there exists a $C^1$ map $\widetilde{H}:\mathbb{R}^d\rightarrow\mathbb{R}^d$ defined in a neighborhood of $z_*$ such that $y_*=\widetilde{H}(z_*)$
and $z_*$ is a solution to the following \ref{NCP} problem:
\begin{equation}
\left\{
\begin{array}
[c]{l}
\widetilde{H}(z)=F(\widetilde{H}(z)-z)-z\smallskip\\
0\leq z\perp \widetilde{H}(z)\geq0.
\end{array}
\right.  \label{eqNCP2}
\end{equation}

\end{itemize}
\end{lemma}
\noindent See \Cref{proof_implicite-fct} for the proof of this lemma.
\begin{remark} \
\begin{enumerate}[(i).]
    \item The condition $F-I$ (respectively, $-(F+I)$) being a $P_0$-map is  actually quite natural since it is exactly the requested assumptions to solve \ref{NCP}, see \Cref{subsec_smoothing_techniques}.
    \item ({\small \text{\ref{NAVE} vs \ref{AVE}}}). At this stage, the reader may have already noticed an analogy with the \ref{LCP} reformulation of \ref{AVE}. Indeed, if $F(x)=Ax-b$, then \ref{NAVE} coincides with \ref{AVE}, and if either $A-I$ or $-(A+I)$ is a $P_{0}$-matrix (which is automatically satisfied if, \textit{e.g.}, the singular values of the matrix $A$ are greater than~$1$), then the functions $H$ and $\widetilde{H}$ are explicitly given by the formulas 
\[
H(y)=My+q \qquad \text{ and } \qquad \widetilde{H}(y)=\widetilde{M}z+\widetilde{q},
\] 
where $M,\widetilde{M}$, $q$ and $\widetilde{q}$ appear in~\eqref{mounir1}. Consequently, in this case, it is possible to solve \ref{AVE} as explained in the introduction.
\end{enumerate}
\end{remark}
\noindent We conclude with a sufficient condition for the uniqueness of solutions to \ref{NAVE}.
\begin{proposition}\label{p0map}
    Let $F: \mathbb{R}^d \to \mathbb{R}^d$. If either $F-I$ or $-(F+I)$ is a $P_0$-map (resp. $P$-map), then for all $x, y \in \mathbb{R}^d$ such that $x \ne y$, there exists an index $i_0 \in \{1, \cdots, d\}$ such that $x^{i_0} \ne y^{i_0}$ and $$|F(x)^{i_0} - F(y)^{i_0}| \ge |x^{i_0} - y^{i_0}|\quad \big(\text{resp. } |F(x)^{i_0} - F(y)^{i_0}| > |x^{i_0} - y^{i_0}|\big).$$
Especially, if either $F-I$ or $-(F+I)$ is $P$, then the problem \ref{NAVE}
admits at most one solution.
    \end{proposition}
\noindent See \Cref{proof_p0map} for the proof of this proposition.
In general, the converse in Proposition~\ref{p0map} fails: for $F(x) = - \frac{1}{2} |x|$ on $\mathbb R$, the equation $F(x) - |x| = 0$ has a unique solution $x_\ast = 0 $ while neither $F - I$ nor $-(F + I)$ is a $P$-map.




\medskip

As already mentioned, even though the functions $H$ and $ \widetilde{H}$ are only implicitly defined in Lemma~\ref{implicite-fct}, we can still solve~\eqref{eqNCP1}--\eqref{eqNCP2} numerically (we shall do so below), whenever it is guaranteed that 
$H$, $\widetilde{H}$ are $P_{0}$-maps. This is the aim of the following lemma, yielding a criterion based on $F$ given at the beginning.

\begin{lemma}[Guaranteeing $P_{0}$-property for $H$, $\widetilde{H}$]\label{guarant-Pmap} \phantom{-}
\begin{itemize}
 \item[$(i)$.] If $F-I$ is a $P_{0}$-map (resp. $P$-map), then so is $H$ given by Lemma~\ref{implicite-fct}.
\item[$(ii)$.] If $-(F+I)$ is a $P_{0}$-map (resp. $P$-map), then so is $\widetilde{H}$  given by Lemma~\ref{implicite-fct}. 
\end{itemize}
\end{lemma}
\noindent See \Cref{proof_guarant-Pmap} for the proof of this lemma. 




\medskip

To derive error estimates for the algorithm, it is necessary to introduce a quantitative version for $P_0$- and $P$-map, see \Cref{merged_estimates}.
To this end, we define the so-called uniform $P$-map with a given modulus.

\begin{definition}
Let $H = (H^1, \cdots, H^d)^T : \mathbb{R}^d \to \mathbb{R}^d$ and $h:\mathbb{R}_{+}\to\mathbb{R}_{+}$ satisfying $h(0)=0<h(t)$ for every $t>0$.
Then, the mapping $H$ is said to be  a \textit{uniform $P$-map with modulus} $h$ if 
    \begin{align*}
        \max\limits_{1 \leq i \leq d} (x^i - y^i)(H(x)^i - H(y)^i) \geq h(\|x - y\|), \quad\text{ for every } x, y \in \mathbb{R}^d.
    \end{align*}

    
\end{definition}

\section{Description of the numerical method}\label{sec_method}
\subsection{Smoothing techniques to solve \ref{NCP}}\label{subsec_smoothing_techniques}
\noindent To solve~\ref{NCP}, we will apply the smoothing approach proposed in \cite{HM2014} and more precisely the non-parametric technique introduced in \cite{OHBA2022}. The overall approach of \cite{HM2014} is based on smoothing functions $\theta$ satisfying the following properties:
\begin{itemize}
 \item the function $ \theta: \R \to (-\infty, 1)$ is smooth, concave and increasing; 
 \item $\theta(t) <0$, for all $t\in (-\infty, 0)$, $\theta(0)=0$ and $\lim\limits_{t \to + \infty} \theta(t)=1$.
\end{itemize}

Two smoothing functions usually used in practice are:
    \begin{itemize}
        \item  the rational function $\theta_1:\R\to(-\infty,1)$ defined by
        \begin{align}\label{theta1}
            \theta_1 = \dfrac{t}{t+1}, \quad\text{for } \, t\ge 0 \quad\text{and}\quad \theta_1 = t,\quad \text{for }\, t\le 0.
        \end{align}
        \item  the exponential function $\theta_2:\R\to(-\infty,1)$ defined by
        \begin{align}\label{theta2}
            \theta_2 = 1-e^{-t},\, \text{ for all } t\in\mathbb{R}.
        \end{align}
    \end{itemize}


\begin{figure}[H]
  \centering
  \begin{subfigure}[b]{0.49\linewidth}
    \includegraphics[width=\linewidth]{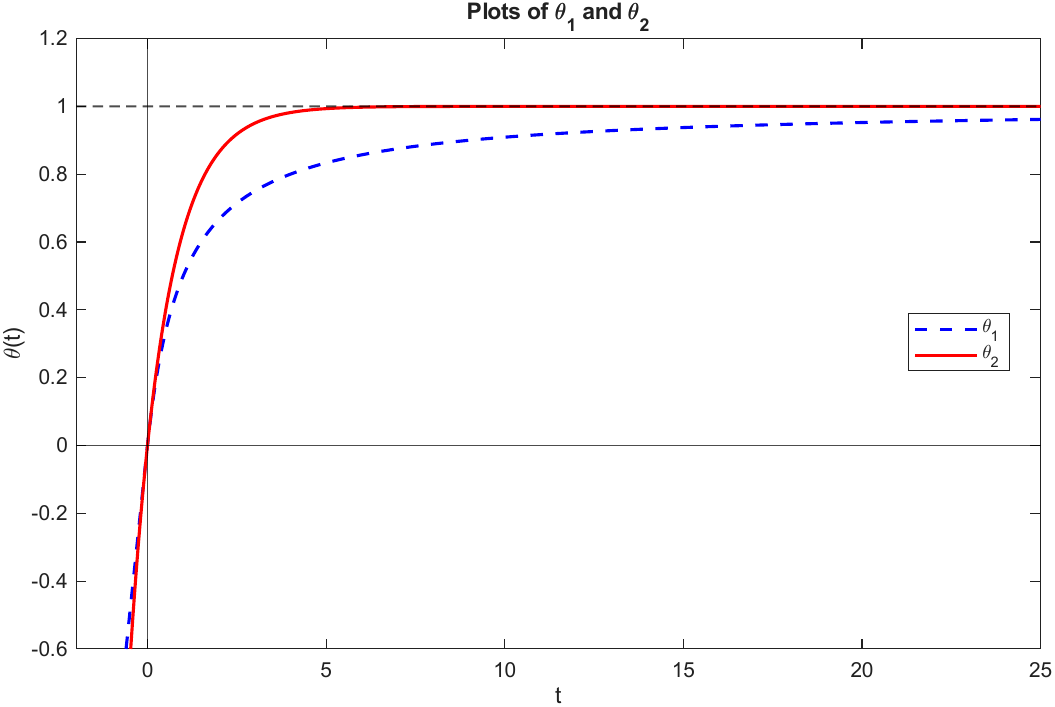}
  \end{subfigure}
  \begin{subfigure}[b]{0.49\linewidth}
    \includegraphics[width=\linewidth]{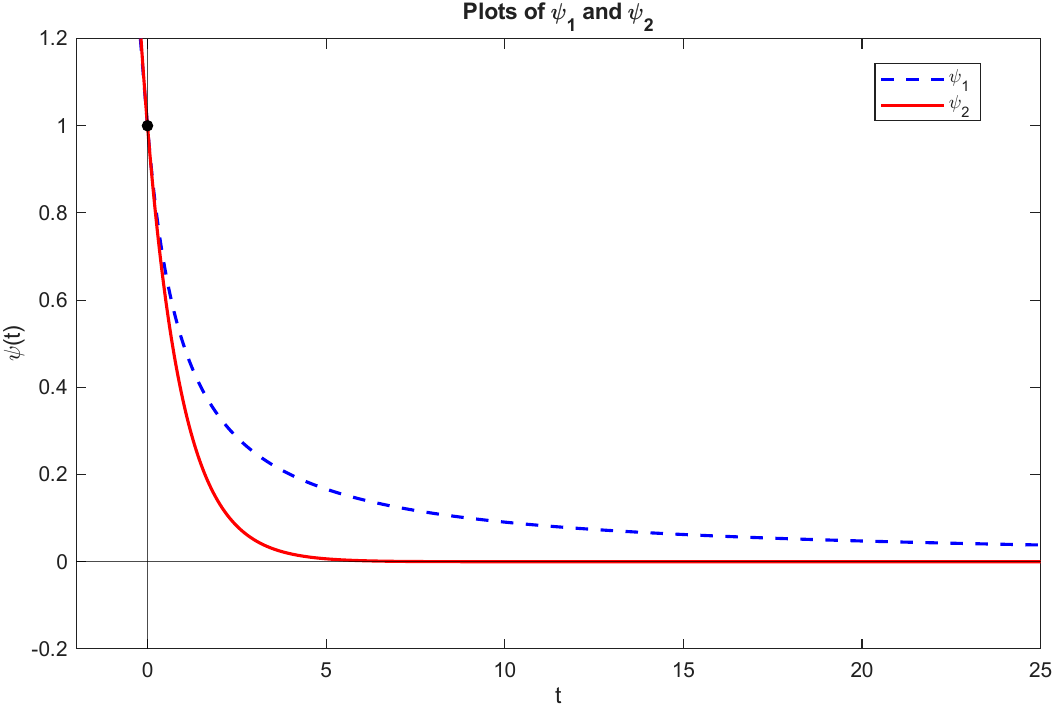}
  \end{subfigure}
  \caption{The graphs of $\theta_i$ and $\psi_i\coloneqq 1-\theta_i$ with $i=1,2$.}
  \label{fig:smoothing-func-plots}
\end{figure}
    
Smoothing functions are used as \textit{certificate of positivity}, that is, they ``detect” whether $t=0$ or $t>0$ holds in a ``continuous way”, in the sense of the following characterization 
$$t>0\quad  \Longleftrightarrow \quad \lim\limits_{r \to 0}\theta\left( \frac{t}{r} \right)=1. $$
The authors in \cite{HM2014} use these functions to regularize the (nonsmooth) \ref{NCP} constraint
\begin{equation}\label{NCP-cstr}
	0 \leq x \perp H(x) \geq 0,
\end{equation}
by means of a sequence of {{smooth}} systems (indexed by $r>0$) of the form 
\begin{equation} \label{smooth-cstr}
 G_{r}(x,H(x))= \left( G_{r}(x,H(x))^1, \cdots ,  G_{r}(x,H(x))^d \right)^T= (0,\cdots ,0)^T,
\end{equation}
where $$ {{ G_{r}(x,H(x))^i:=  r \psi^{-1} \left( \psi \left(\frac{x^{i}}{r}\right)+\psi \left(\frac{H(x)^i}{r}\right) \right)\qquad}}
\text{with}\quad \psi:=1-\theta.$$ Then they eventually take the limit as $r$ tends to~$0$,
proving that solutions to~\eqref{smooth-cstr} converge to a solution of~\eqref{NCP-cstr}, see~\cite[Theorem~4.1]{HM2014}.


\medskip

Several convergence results have been established under the assumption that the problem has at least one solution and $H$ is a $P_0$-map. Although this approach is efficient numerically, it suffers from two drawbacks:
 \begin{itemize}
 \item There is no clear or optimal strategy to drive the parameter $r$ to~$0$.
 \item The following ad hoc technical assumption on the function $\psi$ has been used
to make the algorithm converge
   without rigorous explanation: 
\begin{eqnarray}\label{hyp-techn}
   \text{there exist $a\in (0,1)$ and $R_{a}>0$ such that: \,
$\frac{\psi(t)}{2}\, \geq \, \psi \left(\frac{t}{a} \right)$, \, for all $t \in (R_{a},+\infty)$.}
\end{eqnarray}  
 \end{itemize}
The first drawback has been addressed in~\cite{OHBA2022} by considering a larger system of equations
where $r$ becomes an unknown,
\begin{equation} \label{imp-aug_eq}
\left\{
\begin{array}{l}
G_{r}(x,H(x))=  (0,\cdots ,0)^T,\\[2mm]
\frac{1}{2} \Vert x^{-}\Vert^{2}+\frac{1}{2} \Vert H(x)^{-}\Vert^{2}+ r^{2}+\varepsilon r=0\,.
\end{array}
\right.
\end{equation}
The idea is that,
looking for a solution $(x,r)$ by means of a Newton-like method for some fixed positive parameter $\varepsilon >0$,
will lead to a solution $x$ of~\eqref{NCP-cstr} and $r=0$. Note that the second equation in~\eqref{imp-aug_eq}
ensures that $r$ will quadratically converge to 0.

\medskip

The second drawback will be the subject of the following result which proves that this technical assumption~\eqref{hyp-techn} corresponds to a well-known property.

\begin{theorem}
[asymptotic behavior]\label{thm-mounir} Let $\psi:(0,\infty)\rightarrow
(0,\infty)$ be a smooth convex decreasing function satisfying
\[
\underset{x\rightarrow +\infty}{\lim}\psi(x)=\inf \psi=0.
\]
The following assertions are equivalent:
\begin{enumerate}[(i).]
\item ({\L}ojasiewicz inequality at infinity) There exists $c>0$ such that
\[
\underset{x\rightarrow +\infty}{\lim\inf\,}\,\,\frac{x|\psi^{\prime}(x)|}{\psi(x)}\,\ge\, c>0.
\]
\item There exist $m,n>1$ and $R>0$ such that:
\begin{equation}
\frac{\psi(x)}{m}\geq \psi(nx),\qquad\text{for all }x\in(R,+\infty)
\label{eq:mounir}
\end{equation}

\item For every $m>1$ there exist $n>1$ and $R>0$ such that:
\[
\frac{\psi(x)}{m}\geq \psi(nx),\qquad\text{for all }x\in(R,+\infty).
\]
\end{enumerate}
\end{theorem}
\noindent See \Cref{proof_thm-mounir} for the proof of this theorem.

\medskip

\noindent Notice that the technical assumption~\eqref{hyp-techn} corresponds to (ii). Therefore, the above result shows that it is equivalent to assume that $\psi$ satisfies the {\L}ojasiewicz inequality at infinity (i).
This latter condition is often easily verifiable. In particular, it is straightforward
to check that it holds for $\psi_i\coloneqq 1 - \theta_i$ where $\theta_i$, $i=1,2$ is defined
by~\eqref{theta1} and~\eqref{theta2}.
More generally, the result can be stated for nonsmooth convex functions, and (i) is always satisfied if 
the function $\psi$ is semi-algebraic: indeed, in this case, the corresponding Hardy field (that is, the field of germs of 
real semi-algebraic functions at infinity) has rank one, and consequently, for any non-ultimately zero semi-algebraic function $\psi$ in the single variable $x$, the function $x\mapsto x \psi^{\prime}(x)/\psi(x)$ has a non-zero limit as $x$ goes to infinity (see~\cite[Remark~2.9]{DG2005}). The same argument also applies to the more general case of functions $\psi$ that are definable in some polynomially bounded o-minimal structure (we refer to \cite{Coste2000} for the corresponding definitions). 
This already provides a broad assembly of examples of functions satisfying~(i), together with straightforward criteria to detect easily whether the property holds, based on 
certificates of semialgebricity or o-minimality (see \cite[Theorem~1.13]{Coste2000} \textit{e.g.}).

\begin{remark} Let us draw the reader's attention to the fact that besides what is asserted in \cite[Proposition~2.7]{DG2005}, the assumption of \textit{polynomial boundedness} is essential for the validity of~(i). Indeed, as shown in \cite[Remark~8]{Loi2016}, the convex function $\psi(x)=\left(  \log(1+x)\right)^{-1}$ is definable in the $\log$-$\exp$ o-minimal structure but fails to satisfy~(i).
\end{remark}
\subsection[]{Algorithm and error estimates}\label{subsec_algorithm}


We now describe our algorithm to solve \ref{NCP}. We assume in the following that
case (i) of Lemma~\ref{implicite-fct} (one can proceed in a similar way when case (ii) holds).
We approximate the complementary constraint~\eqref{NCP-cstr} with~\eqref{smooth-cstr}.
Since $H$ is defined implicitly, we introduce a new variable, which plays the role
of $H(y)$. Finally the system to solve reads
\begin{eqnarray}\label{ncp-reform}
\mathbb{H}(\mathbb{X})=(0,0,0)^T,
\end{eqnarray}
where
\begin{eqnarray*}
\mathbb{X}:=(y,z,r)^T
\quad\text{and}\quad
\mathbb{H}(\mathbb{X})
:=
\begin{dcases}
z-F(y-z)+y& \\
G_r(y,z)^i= r \psi^{-1} \left( \psi \Big( \frac{y^{i}}{r} \Big) + \psi \Big( \frac{z^i}{r} \Big) \right),& 1\leq i\leq d,\\
\frac{1}{2} \Vert y^{-}\Vert^{2}+\frac{1}{2} \Vert z^{-}\Vert^{2}+ r^{2}+\varepsilon r.& 
\end{dcases}
\end{eqnarray*}
The following algorithm corresponds to a Newton method with standard Armijo line search.

\begin{table}[H]
\centering
\begin{tabular}{llll}
	\rule{1\linewidth}{0.6mm}  \\
	$ \mathbf{Algorithm }$  \\
	\rule{1\linewidth}{0.6mm}
	\\
	1.~ Chose  $\mathbb{X}_{0} =\left(y_0,z_0,r_{0}\right)$ such that, $y_0>0$, $z_0>0$ and $r_{0}=\left\langle y_{0},z_{0}\right\rangle /d$.\\\hspace*{0.5cm} Set $\tau \in (0,1/2)$, $\varrho \in(0,1)$ and $k=0.$ \\
	2.~ If $ \mathbb{H}\left(\mathbb{X}_{k}\right)=0,$ stop.\\
	3.~ Find a direction $ \mathbf{d}_{k} \in \R^{2d+1}$ such that\vspace{0.3cm} \\
	\vspace{0.3cm}
	\hspace*{5cm}$\mathbb{H}\left(\mathbb{X}_{k}\right)+\nabla_{\mathbb{X}}\mathbb{H}\left(\mathbb{X}_{k}\right)\mathbf{d}_{k}=0.  $\\
	4.~ Choose $\varrho^{j_{k}}\in (0,1),~ $ where  $j_{k} \in \N $ is the smallest integer such that
        \vspace{0.3cm}
        \\
        \vspace{0.3cm}
	\hspace*{5cm}
        $\left\{\begin{array}{l}
        \mathbb{X}_{k}+\varrho^{j_{k}}\mathbf{d}_{k}>0\\[2mm]
        \Theta\left(\mathbb{X}_{k}+\varrho^{j_{k}}\mathbf{d}_{k}\right)
        \leq  \Theta\left(\mathbb{X}_{k}\right)
        + \tau\varrho^{j_{k}}~ \nabla \Theta\left(\mathbb{X}_{k}\right)^{T}\mathbf{d}_{k}.       
        \end{array}\right.$\\  
	5. ~Set $ \mathbb{X}_{k+1}=\mathbb{X}_{k}+\varrho^{j_{k}}\mathbf{d}_{k}~ $ and $ k\gets k+1.~ $Go to step $ 2.$ \\
	\rule{1\linewidth}{0.6mm}
	
\end{tabular}
\end{table}

\noindent  The merit function used in the line search corresponds to the square of the global error:
\begin{equation*}
\Theta(\mathbb{X})=\frac{1}{2}\Vert \mathbb{H}(\mathbb{X})\Vert^{2}.
\end{equation*}
By choosing $y_0>0$ and $z_0 >0$, we are sure that $r_0>0$ and that our algorithm is well-defined. 
 
Convergence results when the technical condition~\eqref{hyp-techn} holds and
for the ideal case, in which the first two equations of the \ref{NCP} reformulation~\eqref{ncp-reform} are assumed to be solved exactly,
can be found in \cite[Proposition 4.1]{HM2014}.
Instead, we establish error estimates for the proposed algorithm in
a more practical case where controlled computational errors are tolerated at each iteration.
The precise relaxed setting is given in the following assumption:

\begin{assumption}\label{assumption1} \ 
\begin{itemize}
\item[(i)] $F : \mathbb{R}^d \to \mathbb{R}^d$ is such that \ref{NAVE} has at least one solution and $F$ is $C^1$ in a neighborhood of the solution. In addition, $F-I$ is a $P_0$-map.
\item[(ii)] $\theta$ is a smoothing function such that $\psi=1-\theta$ satisfies~\eqref{hyp-techn} and $\psi\leq \psi_1$.
\item[(iii)] $\{y_k\}_{k\in\mathbb{N}}$ and $\{z_k\}_{k\in\mathbb{N}}$ are sequences of non-negative solutions of
$$
  \big\|z_k - F(y_k - z_k) + y_k \big\| \le  \e_k^{(1)}, \quad
  G_{r_k}(y_k, z_k) \le  \e_k^{(2)} \quad \text{(componentwise)},
$$
where and $\e_k^{(1)}=O(r_k^\gamma)$  and $\e_k^{(2)} =O(r_k^\sigma)$ with $\gamma, \sigma\geq 1$
are computational errors that we assume to be controlled by $r_k$.
\end{itemize}
\end{assumption}

We assume that the computational errors are controlled by $r_k$, which we know converges quadratically to 0,
at each step of the algorithm.
This is necessary to be able to apply the smoothing method. 
The best case is when $\gamma= \sigma =2$, which allows to obtain
optimal estimates as in the case without any error on computations.

\medskip

A key consequence
of the complementarity error $G_{r_k}(y_k, z_k) \le O(r_k^\sigma)$ with $\sigma\geq 1$
is the strict positivity of $\tau$ defined in
the following lemma, the proof of which is given in~\Cref{proof_lem_tau}.

\begin{lemma}\label{lem_tau}
Under Assumption~\ref{assumption1}, there exists $\tau\in (0,1]$ such that
\begin{eqnarray}\label{minor-psi}
\min_{1\leq i\leq d, \ k\in\N} \{ \psi\left(y_k^i / r_k\right)+ \psi\left(z_k^i / r_k\right)\} \geq \tau.
\end{eqnarray}
\end{lemma}

\begin{lemma}\label{lem:H_global}
   Under Assumption \ref{assumption1}, assume that $F - I$ is a uniform $P$-map with the modulus $h(t) \geq t^\beta$ for some $\beta > 0$.
   Let $H: U \to V$ be as in Lemma~\ref{implicite-fct}-(i).
   Then, the following assertions hold true:
   \begin{itemize}
       \item[(i)] $H$ is also a uniform $P$-map;
       \item[(ii)] if $F: \R^d \to \R^d$ is continuous and $\beta > 1$, then $F + I: \R^d \to \R^d$ is a homeomorphism and $H$ admits a continuous extension to the entire space, which is defined via ${H(y) = y - (F + I)^{-1}(2y)}$.
   \end{itemize}
\end{lemma}
See Subsection~\ref{proof:H_global} for the proof of this lemma.

\medskip

Since the computational errors are controlled by $r_k$,
all the estimates will depend on $r_k$ and on $\tau$ defined in~\eqref{minor-psi}.

\begin{proposition}\label{1&2eq-est-bis}
Assume that Assumption~\ref{assumption1} hold.
Let $\tau \in (0, 1]$ be as in Lemma \ref{lem_tau}. 
Then, it holds for each $i\in \{1,\cdots, d\}$ and $k\in\mathbb{N}$
\begin{enumerate}[(i).]

\item (Complementarity errors) There exists $C\ge 0$ such that
\begin{eqnarray}\label{compl-err-i}  
\min\left\{ y_k^i, z_k^i \right\} \le \frac{2 - \tau}{\tau} r_k,\qquad
0\le y_k^i z_k^i \le \frac{2 - \tau}{\tau} {r^2_k} +\frac{C |\psi'(0)|}{\tau}  \left(y_k^i+z_k^i\right) {r^{\sigma}_k}.
\end{eqnarray}

\item (\ref{NAVE} error) Set $x_k:=y_k-z_k$. There exists $C \geq 0$ such that
\begin{eqnarray*}
&&  \left| F({x_k})^i - |{x^i_k}| \right| \le \left( C + \frac{4 - 2\tau}{\tau} \right) r_k. 
\end{eqnarray*}
\setcounter{saveenum}{\value{enumi}}
\end{enumerate}

Suppose, in addition, that $F: \R^d \to \R^d$ is continuous and $F-I$ is a uniform $P$-map with modulus $h(t)\ge t^\beta$ for some $\beta>1$.
Then, it holds for each $i\in \{1,\cdots, d\}$ and $k\in\mathbb{N}$
\begin{enumerate}[(i).]
\setcounter{enumi}{\value{saveenum}}
\item (Approximation on $H$)
\quad$\displaystyle \left\|z_k-H(y_k)\right\| \le O\big(r_k^{\frac{\gamma}{\beta-1}}\big).$

\item the sequences $\{y_k\}$ and $\{z_k\}$ are bounded.

\item (Complementarity errors) If $\gamma > \beta -1$, then
\begin{eqnarray}\label{complem-iv-1}
  \min\left\{ y_k^i, H(y_k)^i \right\} \le O\left(r_k\right) \text{ and }
 y_k^i H(y_k)^i\le O(r_k^{\tilde{\sigma}})
\text{ with } \tilde{\sigma}:= \min \left\{ 2, \sigma , \frac{\gamma}{\beta-1} \right\}.
\end{eqnarray}

\item (\ref{NAVE} error) If $\gamma > \beta -1$, then
$\big| F(y_k-H(y_k))^i - |(y_k-H(y_k))^i| \big| \le O\left(r_k\right).$

\end{enumerate}

\end{proposition}

\noindent See \Cref{proof_1&2eq-est-bis} for the proof of this proposition.
\begin{remark} \ \\
(i) 
Without any assumption on $F$ guaranteeing that all solutions of \ref{NAVE}
are bounded, we cannot avoid that the estimates depend on $y_k, z_k$ as
in the second estimate of~\eqref{compl-err-i}. It means that we need to
be more accurate in the computations when $y_k, z_k$ are large.
When assuming that $F-I$ is a uniform-$P$ map with modulus $h(t)\ge t^\beta$, $\beta>1$,
we know that the solution of \ref{NAVE} is unique (Proposition~\ref{p0map})
and we obtain that $y_k, z_k$ are bounded (an even better, they converge, see Theorem~\ref{merged_estimates}).
\\
(ii) If we assume the stronger condition $\psi \le \psi_2$ in Assumption~\ref{assumption1}\,(ii)
(i.e., if we use a larger smoothing function $\theta$), then
  the error estimates for \ref{NAVE} can be improved.
  More precisely, the constants $O(1/\tau)$ in the estimates of Proposition~\ref{1&2eq-est-bis}\,(i)-(ii) can be replaced
  by $O(\ln \tau)$.
\end{remark}

\medskip



We now state the convergence for the original \ref{NAVE} problem and the corresponding error estimates.

\begin{theorem}\label{merged_estimates}
Assume that $F$ is continuous, $F-I$ is a uniform $P$-map with modulus $h(t)\ge t^\beta$ for some $\beta>1$ and \Cref{assumption1} holds true.
Suppose that $x_*$ is the solution to~\ref{NAVE} and $y_*$ is the solution to \ref{NCP} given by Lemma~\ref{implicite-fct}.
Then, $y_k\to y_*$,
$z_k\to H(y_*)$ and
$x_k=y_k-z_k$ converges to ${x_*}= y_*-H\left(y_*\right)$ with the corresponding estimates
\begin{align*}
    \|y_*-y_k\| \le O(r_k^\alpha) \quad \text{and} \quad \|x_*-x_k\| \le O(r_k^{\tilde{\sigma}/\beta}),
\end{align*}
where $\alpha= \min \big\{ \frac{\tilde{\sigma}}{\beta} , \frac{\tilde{\sigma}}{2} \big\}$
and $\tilde{\sigma}=  \min \big\{ 2, \sigma , \frac{\gamma}{\beta-1} \big\} $.
\end{theorem}
\noindent See \Cref{proof_merged_estimates} for the proof of this proposition.

\begin{remark}\label{uniform_P_strongly_monotone} \ \\
(i) 
Error estimates for smoothing methods for \ref{NCP} are usually derived under uniform monotonicity or the uniform 
$P$-property with a fixed quadratic modulus, see  e.g \cite{HM2014,AC_2022}. 
Although both include strong monotonicity, they are structurally different.
Theorem~\ref{merged_estimates} is more general and gives explicit bounds once $h$ is specified.
\\
(ii) In the ideal case (no computational errors, $\tau = 1$ in Lemma~\ref{lem_tau}), the second equation of~\eqref{ncp-reform}
is exactly solved at every iteration and sharper bounds can be derived, which coincide with
those in \cite[Proposition 4.1]{HM2014} and \cite[Proposition 4.1]{AC_2022}.
However, exact resolution is unrealistic in practice, hence the value of these estimates lies in their robustness to computational errors. When $\sigma=\gamma=2$ and $\beta\leq 2$, we obtain the optimal estimates as in the case without any error on computations.
\end{remark}

\section{Applications}
\label{sec_applications}
In this section, we show that several problems can be naturally restated as \ref{NAVE} and can be solved efficiently thanks to the above transformation. We present numerical experiments with the specific smoothing functions $\theta_1$ and $\theta_2$.
\smallskip\newline
The numerical experiments are conducted in an ordinary computer. All program codes are written and executed in MATLAB R2024a. \noindent
For convenience, we denote $\alpha e^{\beta} \coloneqq \alpha \times 10^{\beta}$ for all numerical results obtained via MATLAB. 
In \Cref{ss3.1,ss3.3}, we employ a similar stopping criterion for every numerical method, by using a tolerance ~$Tol = 1e^{-10}$ and fixing the maximum number of iterations to~$N_{\max} = 2000$.
Since the \ref{NAVE} problems may have multiple solutions, in the following, the error will be computed by $\text{\textit{Error}} = \|F(x_{\mathrm{approximate}}) - |x_{\mathrm{approximate}}| \|$ in Subsection~\ref{ss3.1} and respectively by $\text{\textit{Error}} = \|F(x_{\mathrm{approximate}}) - |x_{\mathrm{approximate}}| - b\|$ in Subsection~\ref{ss3.3}. 

\subsection{Ridge Regression}\label{ss3.1}

Ridge regression adds to the  loss function $\mathcal{L}(x),$ $x\in \mathbb{R}^{d}$ a penalty term in order to avoid overfitting. Historical development and applications in data science of ridge regression can be found, \textit{e.g.}, in \cite{Hoerl_2020, Hastie_2020}. This penalty term usually consists of adding the squared magnitude of the coefficients. \smallskip\newline
We hereby consider an asymmetric ridge regression of the form:
\begin{eqnarray}\label{ridge} 
\underset{x\in\mathbb{R}^{d}}{\min}\,\,\left\{\mathcal{L}(x)+ \sum_{j = 1}^d \left( \lambda^j \max\{x^j, 0\}^2 + \mu^j \max\{-x^j, 0\}^2 \right)\right\},
\end{eqnarray}
where $\mathcal{L}(x)$ is differentiable and the penalization parameters $\lambda^j$ and $\mu^j$ satisfy $\lambda^j - \mu^j \neq 0$ for all $j \in \{1, \cdots, d\}$. The case $\lambda^j = \mu^j$ for every $j$ corresponds to the classical ridge regression, which will not be considered here. On the other hand, the case $\lambda^j = 0$ for all $j$ and $\mu^j > 0$, corresponds to a penalization of the negativity of the coefficients, promoting solutions with positive coefficients. The necessary condition for optimality reads as follows: $$\nabla \mathcal{L}(x)+ 2 \lambda \max\{x, 0\} - 2\mu \max \{-x,0\}=0,$$
where the two vectors~$\lambda\max \{x,0\}$ and~$\mu\max\{-x, 0\}$ are to be understood componentwise.
Noticing that~$2\max\{x, 0\} = |x| + x$ and~$2 \max \{-x,0\}= |x|-x$, we end up with the following \ref{NAVE} problem
\begin{eqnarray*}
F(x)-|x|=0 \qquad \text{ with}\qquad F(x)=\left(\frac{1}{\mu - \lambda}\right)\nabla \mathcal{L}(x)\,+\,\left(\dfrac{\mu + \lambda}{\mu - \lambda}\right)\,x  \,\,\text{~ (coordinatewise)}
\end{eqnarray*}
Therefore, one can solve the previous problem if either $F-I$ or $-(F+I)$ is a $P_0$-map,
that is,
\begin{align}\label{p0_ridge_regression}
    \text{either }\quad (\mu - \lambda)^{-1}\left( \nabla\mathcal{L} + 2\lambda I \right) \quad\text{ or } \quad - (\mu - \lambda)^{-1}\left(\nabla \mathcal{L} + 2\mu I\right) \quad \text{is a }\,\, P_0\text{-map.}
\end{align}
To illustrate for asymmetric ridge regression, we consider the loss function 
\begin{equation}\label{eq:loss}
\mathcal{L}(x) = \frac{1}{2}\|Ax - b\|^2,\qquad \text{where } \,A \in \R^{m \times d} \,\,\text{ and } \,\, b \in \R^d.
\end{equation}
We show in \Cref{guarantee_sec3.1} that~\eqref{p0_ridge_regression} is indeed true if $\mu>\lambda\ge 0$ or $\lambda>\mu\ge 0$.
\smallskip\newline
We perform numerical experiments, fixing $\lambda^j = \bar\lambda$ and $\mu^j = \bar\mu$ for every $j \in \{1, \cdots, d\}$. These parameters, matrix $A \in \R^{m \times d}$ and vector $b \in \R^d$ were randomly generated with values in $[-5, 5]$. \smallskip\newline
{As shown in Table~\ref{Table.Ridge}, considering the average number of iterations with similar tolerance, the smoothing function $\theta_2$ performs well in all tests, consistently achieving errors below $10^{-17}$    and even zero in two cases - demonstrating both efficiency and stability. In contrast, with $(\bar\lambda, \bar\mu) = (0,10^3)$,  $\theta_1$ requires 2000 iterations to reach $10^{-9}$ error for the case $(m,d) = (10,10)$, and fails to converge for $(m,d) = (20,10)$. On the other hand, while the parameters $\bar\lambda$ and $\bar\mu$ become greater, which can be compared to the ascent of the (classical) ridge parameter, $\theta_2$--smoothing performs within a better tolerance in a small number of iterations.}
\begin{table}[H]
\centering
\captionof{table}{Comparing (asymmetric) ridge regression with different smoothing functions}
\begin{tabular}{c c c c c c c c}
\rule[-1ex]{0pt}{2.5ex}  &  & \multicolumn{2}{c}{Error} & \multicolumn{2}{c}{Iterations} & \multicolumn{2}{c}{Running time $\left(\times e^{-2}(\text{s})\right)$}  \\ 
\cline{3-8} 
\rule[-1ex]{0pt}{2.5ex}  $(\bar\lambda, \bar\mu)$& $(m, d)$  & $\theta_1$ & $\theta_2$ & $\theta_1$ & $\theta_2$ & $\theta_1$ & $\theta_2$\\ 
\Xhline{1pt}
\rule[-1ex]{0pt}{2.5ex}  $(0, 1000)$& $(3, 10)$  & $2.1e^{-11}$ & $7.7e^{-17}$ & $20$ & $28$ & $1.58$ & $1.97$  \\ 
\rule[-1ex]{0pt}{2.5ex} & $(5, 10)$ & $1.4e^{-11}$ & $6.2e^{-17}$ & $21$ & $27$ & $1.54$ & $1.65$  \\ 
\rule[-1ex]{0pt}{2.5ex}  & $(10, 10)$  & $2.3e^{-9}$ & $2.9e^{-17}$ & $2000$ & $30$ & $19.3$ & $2.28$  \\
\rule[-1ex]{0pt}{2.5ex}  & $(20, 10)$  & $6.4e^{-1}$ & $7.9e^{-17}$ & $2000$ & $27$ & $428$ & $1.69$  \\
\Xhline{1pt}
\rule[-1ex]{0pt}{2.5ex}  $(200, 10000)$& $(3, 10)$  & $1.4e^{-11}$ & $0$ & $20$ & $23$ & $1.55$ & $1.68$  \\ 
\rule[-1ex]{0pt}{2.5ex} & $(5, 10)$  & $1.4e^{-11}$ & $8.5e^{-18}$ & $20$ & $23$ & $1.55$ & $1.68$  \\  
\rule[-1ex]{0pt}{2.5ex}  & $(10, 10)$  & $1.4e^{-11}$ & $0$ & $18$ & $22$ & $1.6$ & $1.58$  \\
\rule[-1ex]{0pt}{2.5ex}  & $(20, 10)$  & $5.4e^{-12}$ & $6.9e^{-18}$ & $19$ & $23$ & $1.53$ & $1.67$  \\
\Xhline{1pt}
\end{tabular}
\label{Table.Ridge}
\end{table}

\noindent To end this part, we give a heuristic observation on a sparse optimization problem (see, {\it e.g.}, \cite{T_1996, HTF_2001})
\begin{eqnarray}\label{sparse} 
\underset{x\in\mathbb{R}^{d}}{\min }\,\mathcal{L}(x)+ \lambda \|x\|_1. 
\end{eqnarray}
The first order optimality condition for \eqref{sparse} has the form
	\begin{equation}\label{sparse.condi}
		0	\in \nabla \mathcal{L}(x)\, +\, \lambda \,\partial \|\cdot\|_1(x),
	\end{equation}
where the subdifferential of $\ell_1$ norm can be written explicitly as
	\begin{align*}
		q \in \partial \|\cdot\|_1(x) \text{ if and only if } 
		\begin{dcases}
			~ q_i = \text{sign}(x_i), & \text{ if  } x_i \neq 0, \\
			~ \phantom{tri}|q_i| \leq 1, & \text{ if } x_i = 0.
		\end{dcases}
	\end{align*}
Using the fact that $\alpha \, \text{sign}(\alpha) = |\alpha|$ for every $\alpha \in \R$, the inclusion \eqref{sparse.condi} can be transformed into 
	\begin{align*}
		x \nabla \mathcal{L}(x) + \lambda |x| = 0 \, \text{ (coordinatewise). }
	\end{align*}
The above equation just provides a necessary condition for an optimal solution. In the following, we give a short numerical observation to guarantee its potential utility in sparse optimization. In order to apply results from previous sections, it is necessary to ensure that
	\begin{align*}
		- \dfrac{1}{\lambda} \,x \,\nabla \mathcal{L}(x) - I \qquad\text{ or }\qquad \dfrac{1}{\lambda}\,x \,\nabla \mathcal{L}(x) - I \qquad\text{is a  $P_0$-map.}
	\end{align*}
\noindent In the following figures, we use the same quadratic loss function~\eqref{eq:loss}, where matrix $A$ and vector~$b$ are randomly generated ranging from $-1$ to $1$ and $-0,05$ to $0$, respectively. {\Cref{fig:sparsity_new} illustrates the behavior of each coefficient as the tuning parameter \( \lambda > 0 \) increases. In both cases, we observe that the results obtained using the smoothing function \( \theta_1 \) are consistently better than those obtained with \( \theta_2 \)}
\begin{figure}[H]
  \centering
  \begin{subfigure}[b]{0.39\linewidth}
    \includegraphics[width=\linewidth]{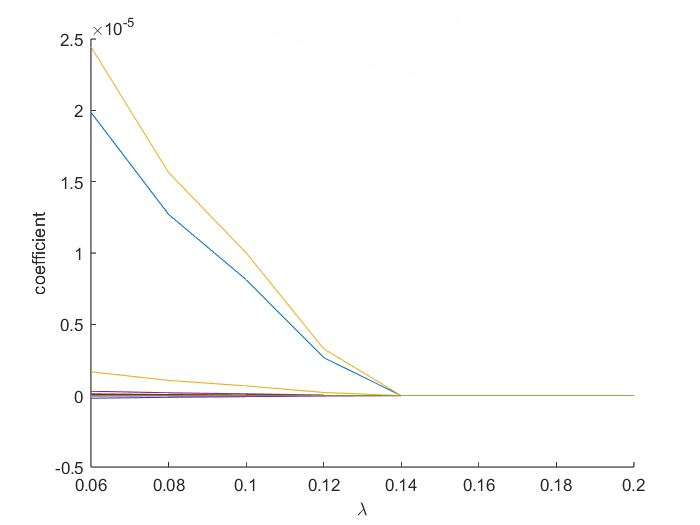}
    \caption{$\theta_1$--smoothing}
  \end{subfigure}
  \begin{subfigure}[b]{0.39\linewidth}
    \includegraphics[width=\linewidth]{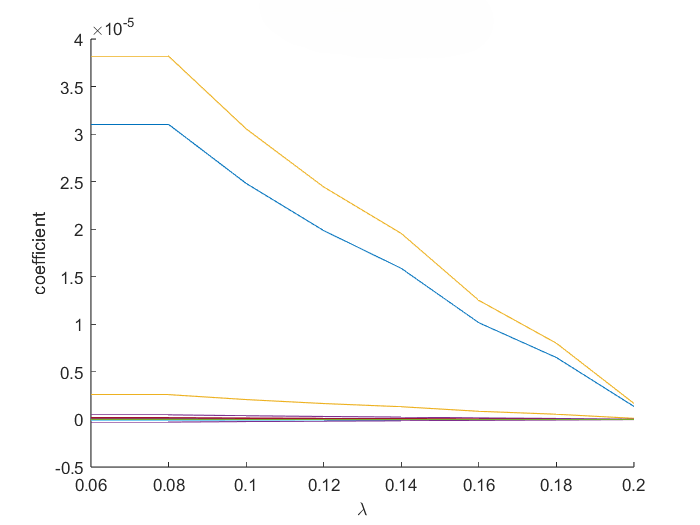}
    \caption{$\theta_2$--smoothing.}
  \end{subfigure}
  \caption{Problem in dimension $m = 1000$ and $d = 2000$.}
  \label{fig:sparsity_new}
\end{figure}

\subsection{Nonlinear ordinary differential equations}
A \ref{NAVE} problem also naturally arises when we deal with a discretization of a nonlinear ordinary differential equation (ODE, for short) involving rough velocity, for example $\dot \gamma(t) = \sqrt{|\gamma(t)|}$ or, more generally, an ODE of the form 
	\[
\Phi(X^{(2k)},X^{(2k-1)},\cdots,\dot{X})=|X|.
	\]
In this subsection, we provide two examples (one being a stiff ODE) to illustrate the effectiveness of smoothing techniques when using finite difference schemes for ODEs.
\begin{example}\normalfont\label{exm.stiff.ODE}
We consider a stiff ODE with initial values as follows
	\begin{equation}\label{stiff.ODE}
		\begin{dcases}
			\ddot{x} + 1001 \dot{x} - 1000|x| = 0, \, t > 0\\
			x(0) = x_0 < 0, \, \dot{x}(0) = 0,
		\end{dcases}
	\end{equation}
whose exact solution is
	\begin{align*}
		x_{\text{exact}}(t) = x_0 \left( - \dfrac{1}{999}e^{-1000t}  + \dfrac{1000}{999} e^{-t} \right) \approx x_0 e^{-t}.
	\end{align*}
Let us consider problem \eqref{stiff.ODE} in time domain $I = [0, T]$. We use a uniform mesh $\pmb{t} = (t_i)$, where $t_i = ih$ for $i \in \{0, \cdots, N \}$ and $h = T/N$, and the approximation solution will be $\pmb{x} = (x_i)$ where $x_i \approx x(t_i)$. For the first and the second derivative, we use the 2nd--order approximation
	\begin{align*}
		\ddot{x}(t_i) \approx \dfrac{x_{i - 2} - 2x_{i - 1} + x_i}{h^2} \text{ and } \dot{x}(t_i) \approx \dfrac{x_{i + 1} - x_{i - 1}}{2h}.
	\end{align*}
Remarkably, at the final time, the first derivative $\dot{x}(t_N)$ will be computed via the 2nd--order backward formula $\dot{x}(t_N) \approx (x_{N - 2} - 4x_{N - 1} + 3x_N)/2h$. Since the initial velocity is zero, using 1st--order backward approximation, we note that $x_{-1} = x_0$. The discretization of \eqref{stiff.ODE} can be written as
	\begin{align*}
		\dfrac{1}{1000}\pmb{A} \pmb{x} + \dfrac{1001}{1000}\pmb{B} \pmb{x} - |\pmb{x}| = \pmb{b},
	\end{align*}
where $\pmb{A}, \pmb{B} \in \R^{N \times N}$ is determined by
	\begin{equation}
		\pmb{A} = \dfrac{1}{h^2}
		\begin{pmatrix}
			1 & 0 & 0 & \cdots & 0 & 0 & 0 \\
			-2 & 1 & 0 & \cdots & 0 & 0 & 0 \\
			1 & - 2 & 1 & \cdots & 0 & 0 & 0 \\
			\cdots & \cdots & \cdots & \cdots & \cdots & \cdots & \cdots \\
			0 & 0 & 0 & \cdots & 1 & 0 & 0 \\
			0 & 0 & 0 & \cdots & -2 & 1 & 0 \\
			0 & 0 & 0 & \cdots & 1 & -2 & 1 
		\end{pmatrix}
		\text{ and } \,\,
		\pmb{B} = \dfrac{1}{2h}
		\begin{pmatrix}
			0 & 1 & 0 & \cdots & 0 & 0 & 0 \\
			-1 & 0 & 1 & \cdots & 0 & 0 & 0 \\
			0 & -1 & 0 & \cdots & 0 & 0 & 0 \\
			\cdots & \cdots & \cdots & \cdots & \cdots & \cdots & \cdots \\
			0 & 0 & 0 & \cdots & 0 & 1 & 0 \\
			0 & 0 & 0 & \cdots & -1 & 0 & 1 \\
			0 & 0 & 0 & \cdots & 1 & -4 & 3 
		\end{pmatrix}
	\end{equation}
Here, vector $\pmb{b} = (b_i) \in \R^N$ is defined by $b_i = 0$ for $i > 2$, $b_1 = x_0(1/(1000h^2) + 1001/(2000h))$ and $b_2 = -x_0/(1000h^2)$.

\noindent In Figure~\ref{fig:stiff}.(a), we approximate the solution of equation \eqref{stiff.ODE} with initial condition $x_0 = -1$ and time interval $I = [0, 4]$. The finite difference scheme was computed with mesh size $h = 0.05$ and the error is $9.22e^{-4}$ when applying $\theta_1$ and $\theta_2$ smoothing functions. To get the convergence rate in Figure~\ref{fig:stiff}.(b), we apply different mesh sizes in the same time interval $I = [0, 1]$ and initial condition $x_0 = -2$. 

\begin{figure}[H]
  \centering
  \begin{subfigure}[b]{0.49\linewidth}
    \includegraphics[width=\linewidth]{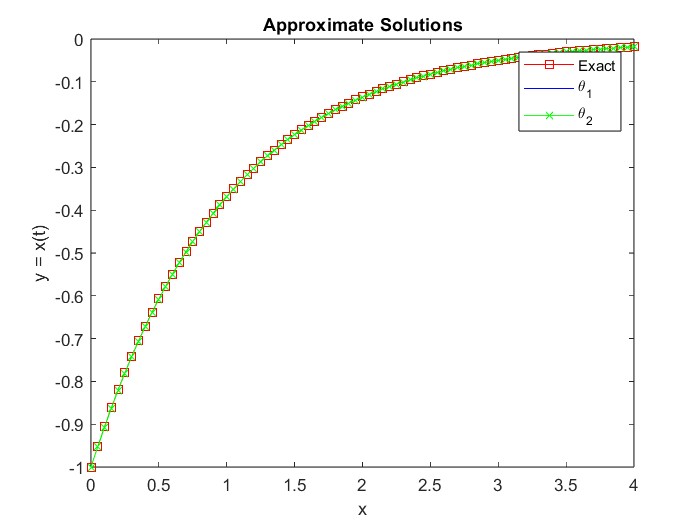}
    \caption{Approximate solutions.}
  \end{subfigure}
  \begin{subfigure}[b]{0.49\linewidth}
    \includegraphics[width=\linewidth]{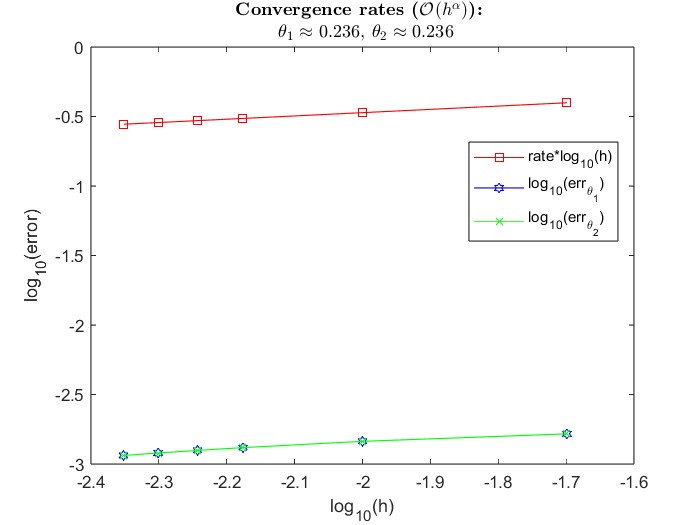}
    \caption{Convergence rate ~$O(h^{1/2})$.}
  \end{subfigure}
  \caption{Solving equation~\ref{stiff.ODE}}
  \label{fig:stiff}
\end{figure}

\noindent Let us now make some comments on the utility of \ref{NAVE} for boundary value problems: we consider a boundary value problem related to \eqref{stiff.ODE} 
	\begin{equation}\label{BVP.ODE}
		\begin{dcases}
			\ddot{x} + 1001\dot{x} - 1000|x| = 0, \, t \in (0, T), \\
			x(0) = x_0 < 0, \, x(T) = y_0 \in \R.
		\end{dcases}
	\end{equation}
In order to illustrate this case, we consider the time interval $I = [0, 2]$ and the exact solution is determined by $x_{\text{exact}}$. Using a similar time mesh as above, the first and second derivatives are approximated as follows
	\begin{align*}
		\ddot{x}(t_i) \approx  \dfrac{x_{i - 1} - 2x_i + x_{i + 1}}{h^2} \text{ and } \dot{x}(t_i) \approx \dfrac{x_i - x_{i - 1}}{h}.
	\end{align*}
Figure~\ref{fig:BVP} shows the convergence rate for the boundary value problem \eqref{BVP.ODE}. The smoothing technique used in this problem presents a better accuracy compared to the above initial value problem, which seems to be natural because of the stiffness of the problem \eqref{stiff.ODE}. It is noteworthy that Figure~\ref{fig:BVP} also depicts an expected convergence rate since we have used a first-order approximation for $\dot{x}$. 
\begin{figure}[H]
	\centering
	\includegraphics[width=0.49\linewidth]{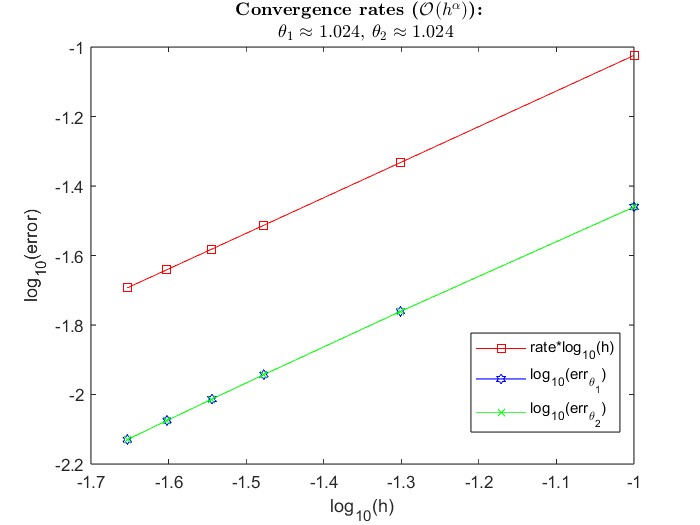}
	\caption{Convergence rate $O(h)$ for a boundary value problem}
  	\label{fig:BVP}
\end{figure}
\end{example}

\begin{example}\normalfont
For a continuous function $f: [0, + \infty) \rightarrow \R$, we consider the ODE
	\begin{equation}\label{actan.ODE}
		\begin{dcases}
			\ddot{x} + \arctan(x) - |x| = f(t), \, t > 0, \\
			x(0) = x_0 \in \R, \, \dot{x}(0) = 0.
		\end{dcases}
	\end{equation}
Using a similar discretization as in Example~\ref{exm.stiff.ODE}, the unknown variable $\pmb{x} \approx x(\pmb{t})$ solves a \ref{NAVE} problem as follows
	\begin{equation}
		\pmb{A}\pmb{x} + \text{arctan}(\pmb{x}) - |\pmb{x}| = \pmb{b},
	\end{equation}
where the matrix $\pmb{A}$ is determined as in Example~\ref{exm.stiff.ODE} and the vector $\pmb{b} \in \R^N$ is defined by $b_i = f(t_i)$ for $i > 2$ and
	\begin{align*}
		b_1 = f(t_1) + \dfrac{x_0}{h^2} \text{ and } b_2 = f(t_2) - \dfrac{x_0}{h^2}.
	\end{align*}
To illustrate for this example, we consider problem \eqref{actan.ODE} with source term $$f(t) = \arctan(\cos(\pi t)) - |\cos(\pi t)| - \pi^2 \cos(\pi t),$$ whose exact solution is $x_{exact}(t) = \cos(\pi t)$. Figure~\ref{fig:arctan}.(a) shows the approximate solution on the time interval $I = [0, 1]$ with mesh size $h = 0.0125$. The error between $\theta_1$ (resp. $\theta_2$) approximation and exact solution is $0.06$ (resp. $0.0725$). Besides, Figure~\ref{fig:arctan}.(b) displays the convergence rate of the combination of finite difference scheme and the $\theta$--smoothing applying for the associated \ref{NAVE} problem. 

\begin{figure}[H]
  \centering
  \begin{subfigure}[b]{0.49\linewidth}
    \includegraphics[width=\linewidth]{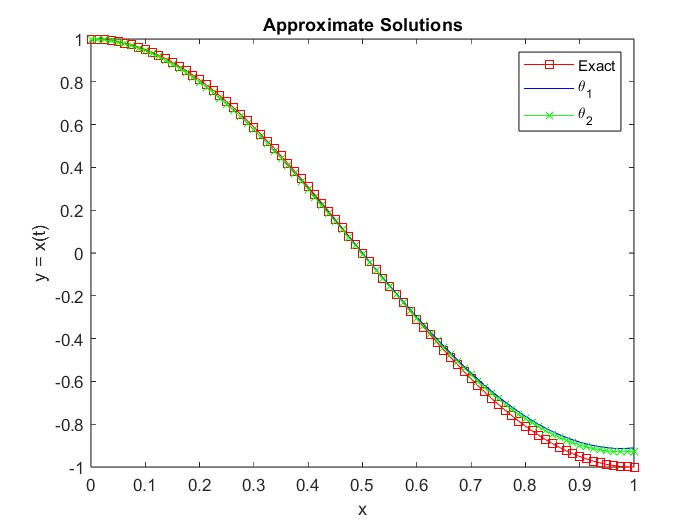}
    \caption{Approximate solutions.}
  \end{subfigure}
  \begin{subfigure}[b]{0.49\linewidth}
    \includegraphics[width=\linewidth]{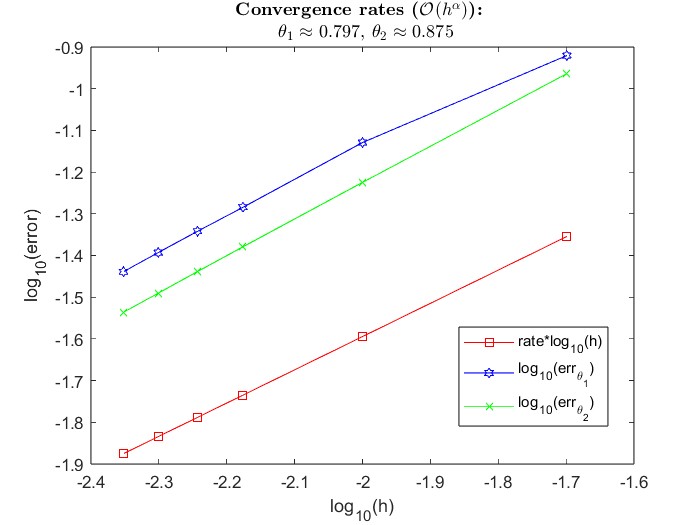}
    \caption{Convergence rate ~$O(h)$.}
  \end{subfigure}
  \caption{Solving equation~\ref{actan.ODE}.}
  \label{fig:arctan}
\end{figure}
\end{example}

\bigskip

\subsection{Comparison of methods for \ref{NAVE}} \label{ss3.3}

\ref{NCP} can be solved with several numerical methods. We give examples to compare the efficiency of the following methods

	\begin{itemize}
		\item The smoothing method developed in \Cref{sec_method}: Newton--like method with smoothing functions $\theta_1$ and $\theta_2$.  
		\item The Soft--Max method SM: approximation by Soft--Max function, in which the main idea is to approximate the complementarity condition via the limit
			\begin{align*}
				\max_{1\le i\le d}  x_i = \lim_{r \searrow 0} r \log \left( \sum_{i = 1}^d e^{x_i/r} \right).
			\end{align*}
This method has been widely used in many optimization problems, for example \cite[Example 1.30]{RW_1998}, \cite{N_2005, L_1992, LTL_2012};
\item The interior point method IP: see for example~\cite{KMNY_1991, PY_1996, I_1997, HMO_2019} for the use of interior point method for complementarity problems.

	\end{itemize}

We solve the system $\tilde{F}(x) - |x| = b$ where $\tilde{F}$ is given  in the following examples.
\begin{example}\normalfont\label{1st_ex}
We consider $\tilde{F}(x) = Ax$, where
	\begin{equation}
		A = \text{tridiag}(-1, 4, -1) \in \R^{d \times d}, \ {x_*} \in \R^d,  \ b = A{x_*} - |{x_*}|.
	\end{equation}
\end{example}
\begin{example}\normalfont\label{2nd_ex}  \cite{KS_1986, AC_2022} $\tilde{F}: \R^3 \rightarrow \R^3$ is defined by
	\begin{align*}
		\tilde{F}(x) := 
			\begin{pmatrix}
			2x_1 - 2 \smallskip \\
			2x_2 + x_2^3 - x_3 + 3 \smallskip \\
			x_2 + 2x_3 + 2 x_3^3 - 3
		\end{pmatrix}.
	\end{align*}
\end{example}
\begin{example}\normalfont\label{3rd_ex}  \cite{KS_1986, AC_2022} $\tilde{F}: \R^4 \rightarrow \R^4$ is defined by
	\begin{align*}
		\tilde{F}(x) := 
		\begin{pmatrix}
			3x_1^2 + x_1 + 2x_1x_2 + 2x_2^2 + x_3 + 3x_4\smallskip \\
			2x_1^2 + x_1 + x_2^2 + x_2 + 10x_3 + 2x_4 \smallskip\\
			3x_1^2 + x_1x_2 + 2x_2^2 + 3x_3 + 9x_4 \smallskip\\
			x_1^2 + 3x_2^2 + 2x_3 + 4x_4
		\end{pmatrix}.
	\end{align*}
\end{example}
\noindent Table~\ref{Table.Merged.FullComparison} compares the four methods: smoothing method with $\theta_1$ and $\theta_2$,  SM method and IP method on the \ref{NAVE} problem associated to \Cref{1st_ex,2nd_ex,3rd_ex}. In the first example, the vector $b$ is randomly generated with values in $[-5, 5]$ and the problem is considered in dimensions $d = 50, 500, 1000$. In \Cref{2nd_ex,3rd_ex}, we respectively consider $b_1 =  (-1, -5, 10)^T$, $ b_2 = (9, -100, 10)^T$, $b_3 = (200, 0, 900)^T$, $b^*_1 = (10, 10, -12, 0)^T$, $b^*_2 = (20, -100, -12, 1)^T$ and $b^*_3 = (200, 10, -5, -5)^T$. We observe that the smoothing method  (especially with $\theta_2$--smoothing function) is the most robust among the considered methods. In connection with convergence speed, the IP method performs much less competitively than the others, while it only reaches $4.1e^{-5}$ after $N = 2000$ iterations. Another point that can be recognized is that the SM method could only solve problems of small size, for example for problems in dimension greater than or equal to $d = 50$, singularities appear after less than $100$ iterations.

\begin{table}[H]
\centering
\captionof{table}{Comparison of Methods in terms of Error, Iterations, and Running time}
\resizebox{\textwidth}{!}{
\begin{tabular}{c c c c c c c c c c c c c c}
\Xhline{2pt}
 & & \multicolumn{4}{c}{Error} & \multicolumn{4}{c}{Iterations} & \multicolumn{4}{c}{Running time ($\times e^{-2}$ (s))} \\
\cline{3-14}
Example & Vector $b$ & $\theta_1$ & $\theta_2$ & SM & IP & $\theta_1$ & $\theta_2$ & SM & IP & $\theta_1$ & $\theta_2$ & SM & IP \\
\midrule
\multirow{3}{*}{\ref{1st_ex}} 
& $d = 50$ & $6.2e^{-12}$ & $4.6e^{-15}$ & NaN & $7.4e^{-3}$ & 30 & 57 & 32 & $2k$ & 4.67 & 5.27 & 6.07 & 104 \\
& $d = 500$ & $4.5e^{-11}$ & $1.2e^{-14}$ & NaN & $1.9e^{-2}$ & 64 & 165 & 93 & $2k$ & 194 & 516 & 406 & 8608 \\
& $d = 1000$ & $1.7e^{-10}$ & $1.7e^{-14}$ & NaN & $4.7e^{-3}$ & 110 & 275 & 154 & $2k$ & 1717 & 4399 & 3522 & 47532 \\
\midrule
\multirow{3}{*}{\ref{2nd_ex}} 
& $b_1$ & $1.2e^{-11}$ & $2e^{-15}$ & $7.4e^{-14}$ & $5.5e^{-2}$ & 15 & 10 & 8 & $2k$ & 1.3 & 1.51 & 1.78 & 12.8 \\
& $b_2$ & $4e^{-12}$ & $2.3e^{-14}$ & $3.5e^{-11}$ & $1.1e^{-2}$ & 24 & 18 & 15 & $2k$ & 1.62 & 1.73 & 1.92 & 13.2 \\
& $b_3$ & $5e^{-12}$ & $1.4e^{-13}$ & $5.4e^{-13}$ & $2e^{2}$ & 212 & 206 & 205 & $2k$ & 2.06 & 2.03 & 3.02 & 13.5 \\
\midrule
\multirow{3}{*}{\ref{3rd_ex}} 
& $b_1^*$ & $1.3e^{-11}$ & $4.2e^{-15}$ & $1.6e^{-14}$ & $7e^{-2}$ & 17 & 13 & 11 & $2k$ & 1.39 & 1.74 & 5.34 & 12.1 \\
& $b_2^*$ & $1.6e^{-12}$ & $1.6e^{-14}$ & $4.3e^{-14}$ & $1.42$ & 28 & 23 & 16 & $2k$ & 1.64 & 1.77 & 1.94 & 72 \\
& $b_3^*$ & $4.5e^{-12}$ & $1.2e^{-13}$ & $2.1e$ & $4.1$ & 52 & 45 & $2k$ & $2k$ & 1.53 & 1.62 & 96.6 & 653 \\
\Xhline{2pt}
\end{tabular}
}
\label{Table.Merged.FullComparison}
\end{table}

\noindent \Cref{fig:time_a,fig:time_b} display the performance time between different methods for \Cref{1st_ex} with the size $n = 20$ and \Cref{2nd_ex}, respectively. We did the observation with $100$ samples, and the vector $b$ is randomly generated with values in $[-10, 10]$. At first sight, the IP method appears to be the slowest one in comparison with the other three methods. As shown in \Cref{fig:time_a}, the $\theta_1$--smoothing method performs the best among all the methods. If we look carefully, in lower dimension as \Cref{2nd_ex}, the SM and $\theta_2$--smoothing methods performs slightly better than $\theta_1$--smoothing method.
\begin{figure}[H]
  \centering
  \begin{subfigure}[b]{0.49\linewidth}
    \includegraphics[width=\linewidth]{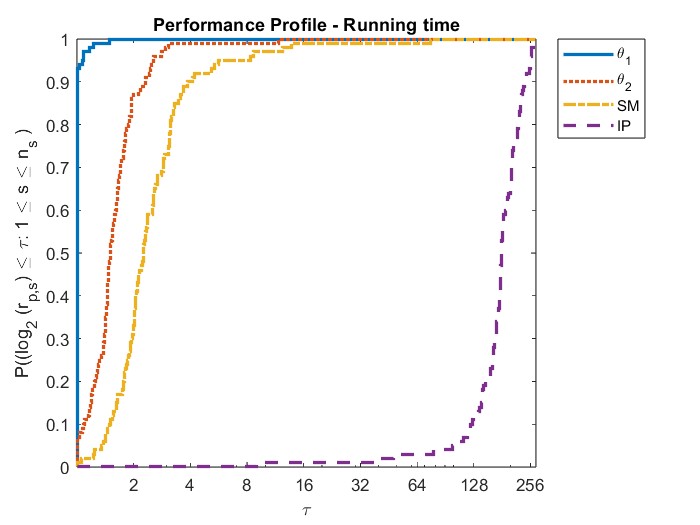}
    \caption{\Cref{1st_ex} in dimension $n = 20$.}
    \label{fig:time_a}
  \end{subfigure}
  \begin{subfigure}[b]{0.49\linewidth}
    \includegraphics[width=\linewidth]{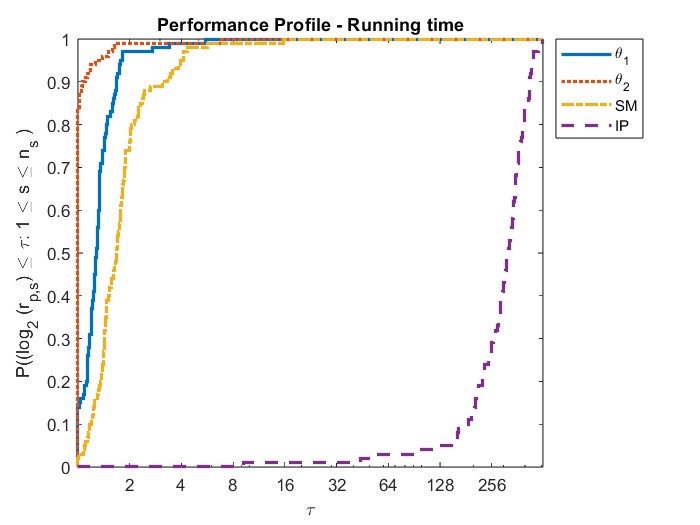}
    \caption{\Cref{2nd_ex}}
    \label{fig:time_b}
  \end{subfigure}
  \caption{Performance time.}
  \label{fig:time}
\end{figure}

\section{Conclusion and discussion}
In this work, we applied smoothing techniques commonly used for Nonlinear Complementarity Problems \ref{NCP} to Nonlinear Absolute Value Equations \ref{NAVE}. We first demonstrated that \ref{NAVE} can be reformulated as an \ref{NCP} with an implicit mapping, allowing us to leverage established \ref{NCP} methods. We subsequently showed that \ref{NAVE} can be effectively solved under weak, direct assumptions regarding the \ref{NAVE} function. Furthermore, we derived algorithm error estimates under the uniform P-property for both \ref{NAVE} and \ref{NCP}, with the latter providing more general results while relaxing standard regularity assumptions. Additionally, we clarified a technical assumption common in smoothing approaches by proving its equivalence to a Łojasiewicz inequality, thus establishing its broader applicability. Finally, through illustrative examples and numerical verification, we confirmed the effectiveness and potential of the proposed approach.

\medskip 

Future work will focus on establishing sharper error bounds and studying the complexity of the method, particularly in favorable scenarios, such as when the implicit mapping in the \ref{NCP} formulation is monotone or uniform-$P$.

\bigskip

\noindent\textbf{Acknowledgement.}  This work was initiated during a research stay of Aris Daniilidis and Tr\'i Minh L\^e to INSA Rennes. These authors thank their hosts for
hospitality. The first author acknowledges support from the Austrian Science
Fund (FWF, P--36344-N).

\smallskip
\noindent\rule{6cm}{1.5pt}


\noindent\rule{6cm}{2pt}
\bigskip


\noindent Aris Daniilidis, Tr\'i Minh L\^e 
\smallskip

\noindent Institut f\"{u}r Stochastik und Wirtschaftsmathematik, VADOR E105-04
\newline TU Wien, Wiedner Hauptstra{\ss }e 8, A-1040 Wien\smallskip
\newline\noindent E-mail: \{\texttt{aris.daniilidis, minh.le\}@tuwien.ac.at}
\newline\noindent\texttt{https://www.arisdaniilidis.at/}
\newline\noindent\texttt{https://sites.google.com/view/tri-minh-le/}

\smallskip

\noindent Research supported by the grants: \smallskip\newline Austrian
Science Fund (FWF P-36344N) (Austria)\newline

\vspace{0.2cm}

\noindent Mounir Haddou, Olivier Ley \& Phi Hoang Tran \smallskip

\noindent Univ Rennes, INSA, CNRS, IRMAR - UMR 6625, F-35000 Rennes, France
\smallskip\newline\noindent E-mail: \texttt{\{mounir.haddou,
olivier.ley, Hoang-phi.Tran\}@insa-rennes.fr} \newline\noindent\texttt{http://\{haddou,
ley\}.perso.math.cnrs.fr/} \smallskip

\noindent Research supported by the Centre Henri Lebesgue ANR-11-LABX-0020-01.
\appendix
\section{Appendix: Proof of the main results}
\subsection[]{Proof of \Cref{lem-P0-invert}}\label{proof_lem-P0-invert}
Let $A$ be a $P_0$-matrix. Then for every diagonal matrix $\Delta_2$ with nonnegative entries, the matrix $\Delta_2 A$ is also $P_0$, while for every diagonal matrix $\Delta_{1}$ with strictly positive entries, the matrix $\Delta_1 + \Delta_2 A$ is a $P$-matrix, therefore, in particular, it is invertible. \smallskip\newline
Conversely, let us assume that $A$ is not a $P_0$-matrix. Then there exists $v\in\R^d$, $v\not= 0$ such that
\begin{eqnarray}\label{contra154}
(Av)^i v^i <0, \qquad \text{for every $i\in\{1,\cdots ,d\}$}.
\end{eqnarray}
Let $\Delta_1=\text{diag}(\delta_1^1, \cdots , \delta_1^d)$ and $\Delta_2=\text{diag}(1, \cdots , 1) =I$. Then
$$(\Delta_1 +\Delta_2 A)v = (\delta_1^1 v^1 + (Av)^1 , \cdots , \delta_1^d v^d + (Av)^d)^T,$$
and by setting, for every $i$, $\delta_1^i := -(Av)^i/v^i$ (which is well-defined and strictly positive thanks to~\eqref{contra154}) we deduce that $(\Delta_1 +\Delta_2 A)v=0$ and therefore $\Delta_1 +\Delta_2 A$ is not invertible. \hfill$\Box$
\subsection[]{Proof of \Cref{implicite-fct}}\label{proof_implicite-fct}
 We consider a $C^1$ map $\mathcal{F}: \mathbb{R}^{2d}\to \mathbb{R}^d$ defined
by $$\mathcal{F}(y,z)=F(y-z)-(y+z).$$ We are going to apply the implicit function theorem at the point $(y_*,z_*)\in \mathbb{R}^{2d}$. Notice that 
\begin{align*}
    \left(\nabla_y\mathcal{F}(y_*,z_*), \nabla_z\mathcal{F}(y_*,z_*)\right)=
\left(\nabla F(y_*-z_*)-I , -\nabla F(y_*-z_*)-I\right).
\end{align*}
If $F-I$ is a $P_0$-map, then $\nabla F({x_*})-I$ is a $P_0$-matrix by Lemma~\ref{lem-map-mat}.
Applying Lemma~\ref{lem-P0-invert}, we obtain that $2I+\nabla F({x_*})-I= \nabla F({x_*})+I$
is invertible. Therefore $\nabla_z\mathcal{F}(y_*,z_*)$ is invertible. By Implicit Function Theorem, there exists a $C^1$ map $H:U\to V$ where $U\ni y_*$ and $V\ni z_*$ are two open subsets of $\mathbb{R}^d$, such that for each $y\in U$, there is a unique $z=H(y) \in V$ satisfying $\mathcal{F}(y,z) =0$. Hence $\mathcal{F}(y_*,z_*) =0$ with $z_* = H(y_*)$. Rewriting this equality and combining it with the fact that $0\le y_* \perp z_* \ge 0$, we conclude that $y_*$ is a solution to NCP~\eqref{eqNCP1}.

\medskip

Similarly, if $-(F+I)$ is a $P_0$-map, then $2I-(\nabla F(z_*)+I)= I- \nabla F(z_*)= -\nabla_y\mathcal{F}(y_*,z_*)$
is invertible and we obtain a $C^1$ map $\widetilde{H}$ such that NCP~\eqref{eqNCP2} holds.  \hfill$\Box$
\subsection[]{Proof of \Cref{p0map}}\label{proof_p0map}
Let $x, y \in \mathbb{R}^d$ be such that $x \ne y$. 
Suppose that $F - I$ is a $P_0$-map. 
Then, by definition, there exists an index $i_0 \in \left\{1, \cdots, d\right\}$ such that $x^{i_0} \ne y^{i_0}$ and
\begin{equation}\label{FI}
\Big( \left(F(x)^{i_0} - x^{i_0}\right) - \left(F(y)^{i_0} - y^{i_0}\right)\Big)\left(x^{i_0} - y^{i_0}\right) \ge 0.
\end{equation}
Rearranging terms in \eqref{FI}, a direct computation yields
\begin{align*}
    | F(x)^{i_0} - F(y)^{i_0}| |x^{i_0} - y^{i_0}| \geq (F(x)^{i_0} - F(y)^{i_0})(x^{i_0} - y^{i_0}) \geq (x^{i_0} - y^{i_0})^2.
\end{align*}
Therefore, we obtain
\begin{equation}\label{FI2}
    |F(x)^{i_0} - F(y)^{i_0}| \ge |x^{i_0} - y^{i_0}|.
\end{equation}
The case where $-(F+I)$ is a $P_0$-map is analogous. 
Replacing $P_0$--property with $P$--property yields a version of~\eqref{FI2} in which the inequality becomes strict.

\medskip

Now assume that either $F-I$ or $-(F+I)$ is a $P$-map and the equation $F(x)=|x|$ has two different solutions $y,z\in\mathbb{R}^d$. 
Following the above observation, there exists $i_0 \in \{1, \cdots, d\}$ such that 
$$|F(y)^{i_0} - F(z)^{i_0}| > |y^{i_0} - z^{i_0}|.$$ 
Since $F(y) = \left|y\right|$ and $F(z) = \left|z\right|$, we infer that
    $$|F(y)^{i_0} - F(z)^{i_0}| > |y^{i_0} - z^{i_0}| \ge \big| |y^{i_0}| - |z^{i_0}| \big| = |F(y)^{i_0} - F(z)^{i_0}|,$$
which is a contradiction.
Therefore,  the equation $F(x) = \left|x\right|$ has at most one solution.\hfill$\Box$
\subsection[]{Proof of \Cref{guarant-Pmap}}\label{proof_guarant-Pmap}
We only prove (i) and the proof of (ii) is analogous.
Let $H: U \to V$ be the map defined as in Lemma~\ref{implicite-fct} - (i).
Fix $y, z \in U$ and set $u := y - H(y)$ and $v := z - H(z)$.
It is straightforward to check that $2H(y) = (F - I)(u)$ and $2H(z) = (F - I)(v)$.
Therefore, for any fixed $i \in \{1, \cdots, d \}$, we have
\begin{equation}\label{Hi_identity}
\begin{split}
& \hphantom{ = \ } \big( 2H(y)^i-2H(z)^i \big)( y^i - z^i)\\ 
&= \big( (F-I)(u)^i - (F-I)(v)^i \big) \big(u^i + H(y)^i - v^i - H(z)^i \big)\\
      &= \big( (F-I)(u)^i - (F-I)(v)^i \big) (u^i - v^i) + \big( (F-I)(u)^i-(F-I)(v)^i \big) \big( H(y)^i - H(z)^i\big)\\
      &= \big( (F-I)(u)^i -(F-I)(v)^i \big) (u^i - v^i )+2\big( H(y)^i -H(z)^i \big)^2,
\end{split}
\end{equation}
Consider the case $y \neq z$.
Since $F-I$ is a $P_{0}$-map, there exists a coordinate
$j=j(u, v) \in \{1,\cdots,d\}$ such that $u^j \ne v^j$ and
\begin{equation}\label{FI_le_0}
\big( (F-I)(u)^j - (F-I)(v)^j \big) \big( u^j- v^j \big) \geq 0.
\end{equation}
The above observation implies that
\begin{equation}\label{ineq386}
\begin{split}
 & ~ 2 (H(y)^j - H(z)^j)(y^j - z^j) \\
 = & ~ \big( (F-I)(u)^i -(F-I)(v)^i \big) (u^i - v^i )+2\big( H(y)^i -H(z)^i \big)^2 \geq 0.
 \end{split}
\end{equation}
To conclude that $H$ is a $P_0$-map, it remains to prove that $y^j \ne z^j$.
Assume by contradiction that $y^j = z^j$. Since
$y^j-H(y)^j = u^j \ne v^j = z^j- H(z)^j$, we infer $H(y)^j \ne H(z)^j$,
which contradicts~\eqref{ineq386}.
Lemma~\ref{guarant-Pmap} is proven.
\hfill$\Box$

\subsection[]{Proof of \Cref{thm-mounir}}\label{proof_thm-mounir}

{(i)$\Longrightarrow$(ii). Suppose that (i) holds.
{Then there exists $R>0$ such that
\[
\inf_{x>R} \frac{x|\psi'(x)|}{\psi(x)} \ge \frac{c}{2}.
\]
}
Hence, using the fact that  $\psi$ is differentiable and decreasing, we obtain
\begin{equation}\label{abs_psi}
    -\frac{\psi'(x)}{\psi(x)} = \frac{|\psi'(x)|}{\psi(x)} \ge \frac{c}{2x}, \text{ for every } x \in (R, +\infty).
\end{equation}
Fix $m > 1$. 
Let $n > 1$ be such that $n^{c/2} \geq m$.
Integrating \eqref{abs_psi} over $(x, nx)$, we get
\begin{align*}
  \ln \psi(x) - \ln \psi(nx) = \int_{x}^{nx} -\frac{\psi'(t)}{\psi(t)} \, dt \ge \int_{x}^{nx} \frac{c}{2t} \, dt
  = \frac{c}{2} \ln (nx) - \frac{c}{2} \ln x, \quad \text{ for every } x > R.
\end{align*}
This implies that $\ln (\psi(x)/\psi(nx)) \geq (c \ln n)/2$. 
Therefore, we have
\begin{align*}
    \frac{\psi(x)}{\psi(nx)} \ge  n^{c/2} \quad \text{ and equivalently } \quad \frac{\psi(x)}{n^{c/2}} \ge  \psi(nx).
\end{align*}
Using the fact that $n^{c/2} \ge m$ and $\psi$ is positive, we deduce
\[
\frac{\psi(x)}{m} \ge \frac{\psi(x)}{n^{c/2}} \ge  \psi(nx), \quad \text{for all } x \in (R, +\infty).
\]

(ii)$\Longrightarrow$(iii). Assume that \eqref{eq:mounir} holds for some
$m_{0},n_{0} >1$ and $R_{0}>0,$ that is, for all $x>R_{0}$ we have
$\psi(x)\,\geq \,m_{0}\,\psi(n_{0}x).$ Then since $n_{0}x>x>R_0$ we also have:
\[
\psi(n_{0}x)\,\geq\, m_{0}\,\psi(n_{0}^{2}x)\quad\text{yielding}\quad \psi(x)\,\geq \,m_{0}
^{2}\,\psi(n_{0}^{2}x).
\]
We conclude that \eqref{eq:mounir} also holds for $m_{1}=m_{0}^{2}$ (under the
choice of $n_{1}=n_{0}^{2}$ and $R_{1}=R_{0}$). Repeating this argument we
deduce that \eqref{eq:mounir} holds for all $m_{k}=m_{0}^{k},$ $k\geq1$
(taking $n_{k}=n_{0}^{k}$ and $R_{k}=R_{0}$). Since $m_{k}\rightarrow\infty,$
in order to establish (iii) it is sufficient to observe that if
\eqref{eq:mounir} holds for some $\bar{m}>1$ (together with some $\bar{n}>1$
and $\bar{R}>0$) then it also holds for all $m\in(1,\bar{m}]$, since
\[
\frac{\psi(x)}{m}\geq\frac{\psi(x)}{\bar{m}}.
\]

(iii)$\Longrightarrow$(i). Fix $m>1,$ $n>1$ and $R>0$ such that
\eqref{eq:mounir} holds and set
\[
c\,:=\,\left(  \frac{m-1}{m}\right)  \left(  \frac{1}{n-1}\right) \, >\,0.
\]
Using convexity of $\psi$ and~\eqref{eq:mounir}, we deduce that for all $x>R$ we have:
\[
\left\vert \psi^{\prime}(x)\right\vert = -  \psi^{\prime}(x)
\,\geq\,\frac{\psi(x)-\psi(nx)}{nx-x}\,\Longrightarrow\,
\frac{x|\psi^{\prime}(x)|}{\psi(x)}\,\geq\,\left(  \frac{1}{n-1}\right)  \left(
1-\frac{\psi(nx)}{\psi(x)}\right) \, \geq \, c.
\]
This establishes (i) and finishes the proof.\hfill$\Box$
\subsection[]{Proof of \Cref{lem_tau}}\label{proof_lem_tau}
Assume that
\begin{eqnarray}\label{lem111}
&& G_{r_k}(y_k, z_k)^i= r_k \psi^{-1} \left( \psi\left(y_k^i / r_k\right)+ \psi\left(z_k^i / r_k\right)\right)\leq \e_k^{(2)}.
\end{eqnarray}
Then, dividing by $r_k$ and applying the decreasing function $\psi$ to the inequality, we obtain
\begin{eqnarray}\label{lem222}
&&  \psi\left(y_k^i / r_k\right)+ \psi\left(z_k^i / r_k\right)\geq \psi(\e_k^{(2)}/r_k ).
\end{eqnarray}
Using the fact that $\e_k^{(2)}\leq Cr_k$ for some $C\geq 0$, we obtain the desired result with $\tau = \psi(C) > 0$.
\hfill$\Box$

\subsection[]{Proof of Lemma~\ref{lem:H_global}}\label{proof:H_global}
(i)
Fix $x,y\in U$ and set $u := x - H(x)$ and $v := y-H(y)$. 
It follows from \eqref{Hi_identity} that
\begin{equation}\label{maxH-01}
\begin{split}
   2\max_{1\le i\le d} \left( H(x)^i - H(y)^i \right) ( x^i - y^i) 
    \ge & ~ \max_{1\le i\le d}\left( (F-I)(u)^i-(F-I)(v)^i\right) ( u^i - v^i)  \\
    \ge & ~  h(\|u-v\|)\ge \|u-v\|^\beta.
\end{split}
\end{equation}
Further, a direct computation yields
\begin{equation}\label{maxH-02}
    \begin{split}
    \left\| u -v \right\|^2 = \, & \left\| x-y \right\|^2 + \left\| H(x)-H(y) \right\|^2 - 2\left\langle x-y,H(x)-H(y)\right\rangle \\
    \geq \, & \left\| x-y \right\|^2 + \left\| H(x)-H(y) \right\|^2 - 2 d  \max\limits_{1\le i\le d} \left( H(x)^i - H(y)^i \right)\left( x^i - y^i \right).
    \end{split}
\end{equation}
Combining \eqref{maxH-01} and \eqref{maxH-02}, we arrive at
\begin{align*}
   (2\lambda)^{2/\beta} + 2d\lambda  \ge \left\| x-y \right\|^2,
\end{align*}
where $\lambda :=\max_{1\le i\le d} \left( H(x)^i - H(y)^i \right)\left( x^i - y^i\right)$. Since $\omega(t)=(2t)^{2/\beta} + 2dt$ is strictly increasing on ${(0,+\infty)}$, the above inequality yields 
\begin{align*}
    \max_{1\le i\le d} \left( H(x)^i - H(y)^i \right)\left( x^i - y^i\right) \ge \omega^{-1}(\|x-y\|^2).
\end{align*}
This proves that $H$ is a uniform $P$-map with the modulus $\widetilde\omega(t) = \omega^{-1}(t^2)$, where $\widetilde\omega(0)\coloneqq 0$.

\medskip

(ii) 
Assume that $F$ is continuous in $\R^d$ and $\beta>1$.
For the sake of brevity, denote $G := F + I$.
We split the proof into two steps.

\smallskip

\textbf{Step 1:} \textit{$G$ is injective and coercive.}
For any $x,y\in \R^d$ with $x\ne y$, the uniform $P$-property of $F-I$ gives
    \begin{equation}\label{G-Esti}
    \begin{split}
        & \, \max\limits_{1\le i \le d}\left(G(x)^i-G(y)^i\right)(x^i-y^i) \\
        = & \, \max\limits_{1\le i \le d}\left\{\left((F-I)(x)^i-(F-I)(y)^i\right)(x^i-y^i) + 2(x^i-y^i)^2\right\} \\
        \ge & \, \max\limits_{1\le i \le d}\left((F-I)(x)^i-(F-I)(y)^i\right)(x^i-y^i) \ge h(\|x-y\|).
    \end{split}
    \end{equation}
Hence, $G$ is also a uniform $P$-map with modulus $h$. 
Moreover, there exists $i_0\in \{1,\cdots,d\}$ such that 
    \begin{align*}
        \left(G(x)^{i_0}-G(y)^{i_0}\right)(x^{i_0}-y^{i_0}) \ge h(\|x-y\|)>0.
    \end{align*}
We infer that $G(x)\ne G(y)$ and so $G$ is injective.

\medskip

Now, applying \eqref{G-Esti} to the case $x \neq 0$ and using the fact that $h(t) \geq t^\beta$, we obtain that
\begin{align*}
    \| G(x) - G(0) \| \| x \| \geq \max\limits_{1\le i \le d}\left(G(x)^i-G(0)^i\right)(x^i - 0) \geq h(\|x\|) \geq \|x\|^\beta.
\end{align*}
Therefore,
\begin{align*}
    \| G(x) \| \geq \| x \|^{\beta - 1} - \| G(0) \|, \quad \text{ for every $x \neq 0$},
\end{align*}
which implies that $G$ is coercive since $\beta > 1$.

\smallskip

\textbf{Step 2:} \textit{$G$ is a homeomorphism.}
Note that $G$ is continuous and injective.
By the invariance of domain theorem, $G(\R^d)$ is open and $G:\R^d\to G(\R^d)$ is a homeomorphism. 

\medskip

We prove that $ G(\mathbb{R}^d)$ is closed. 
Take a sequence $ \{ z_k \} \subset G(\mathbb{R}^d) $ such that $z_k \to z$ as $k \to \infty$.
We aim to show $z \in G(\mathbb{R}^d)$.
For each $k \in \mathbb N$, choose $x_k \in \mathbb{R}^d$ such that $ G(x_k)=z_k.$
Notice that $\{ z_k \}$ is bounded.
Thanks to the coercivity of $G$, we infer that $\{ x_k \}$ is also bounded and so, up to a subsequence, $x_k \xrightarrow{k \to \infty} x$ for some $x \in \R^d$.
By the continuity of $G$, we obtain that $G(x) = z$.
Therefore, \(G(\mathbb{R}^d)\) is closed.

\medskip

\textbf{Conclusion}.
Since the nonempty set $G(\R^d)$ is  both open and closed, we have $G(\R^d) = \R^d$.
Therefore, for any $y \in \R^d$, there exists a unique $z \in \R^d$ such that
\[
    z = y - G^{-1}(2y) = y - (F + I)^{-1}(2y).
\]
With this choice $z = z(y)$, we obtain
\[
    F(y - z) - (y + z) = 0.
\]
Moreover, by the definition, we have
\[
    F(y - H(y)) - (y + H(y)) = 0 \quad \text{ for every $y \in U$.}
\]
The uniqueness of $z = z(y)$ then implies that $H(y) = y - (F + I)^{-1}(2y)$ for every $y \in U$.
Since $(F + I)^{-1}$ is continuous in $\R^d$, $H$ admits a continuous extension defined as follows
\begin{align*}
    H:\R^d &\to \R^d\\
     y &\mapsto H(y) := y-(F+I)^{-1}(2y),
\end{align*}
which completes the proof.
\hfill $\square$
\subsection[]{Proof of \Cref{1&2eq-est-bis}}\label{proof_1&2eq-est-bis}
(i) Fix $k\in\mathbb{N}$ and $1\leq i\leq d$.
Recall that $\psi(t)\leq \psi_1(t)=\frac{1}{1+t}$.
It follows from Lemma~\ref{lem_tau} that
\begin{eqnarray}\label{ineq287}
  &&  \tau\leq \psi(\e_k^{(2)}/r_k ) \leq \psi\left(y_k^i/r_k\right)+\psi\left(z_k^i/r_k\right)
  \leq \frac{1}{1+ y_k^i/r_k} +\frac{1}{1+ z_k^i/r_k}.
\end{eqnarray}
Set for simplicity $\tilde{y} = y_k^i/r_k$ and $\tilde{z} = z_k^i/r_k$.
Applying the inequality $\frac{1}{1+ \tilde{y}} +\frac{1}{1+ \tilde{z}}\leq \frac{2}{1+\min\{ \tilde{y}, \tilde{z}\}}$
to \eqref{ineq287}, we obtain the first estimate in~\eqref{compl-err-i}.
Writing
\[
\frac{1}{1+ \tilde{y}} +\frac{1}{1+ \tilde{z}}
=\frac{2+ \tilde{y}+ \tilde{z}}{1+ \tilde{y} + \tilde{z} + \tilde{y}\tilde{z}},
\]
we get
\begin{eqnarray}\label{betk-1}
  y_k^i z_k^i \le \frac{2 - \beta_k}{\beta_k} {r^2_k} + \frac{1 - \beta_k}{\beta_k} \left(y_k^i+z_k^i\right) r_k
   \le \frac{2 - \tau}{\tau} {r^2_k} + \frac{1 - \beta_k}{\tau}\left(y_k^i+z_k^i\right) r_k,
\end{eqnarray}
where $\beta_k:=  \psi(\e_k^{(2)}/r_k )\geq \tau$.
Further, we have
\begin{eqnarray}\label{betk-2}
&& (1-\beta_k)\leq \frac{1-\psi(\e_k^{(2)}/r_k)}{\e_k^{(2)}/r_k}C r_k^{\sigma -1}
\leq C (-\psi'(0)) r_k^{\sigma -1},
\end{eqnarray}
where the last inequality comes from the convexity of $\psi$.
Plugging this estimate into~\eqref{betk-1} yields  the second estimate in~\eqref{compl-err-i}.

\medskip

(ii) Notice that $|{x^i_k}|  = y_k^i + z_k^i - 2\min\left\{y_k^i,z_k^i\right\}$.
Then, we have
\begin{align}
    \big| F\left({x_k}\right)^i - |{x^i_k}| \big| &= \big| F\left(y_k - z_k\right)^i - \big(y_k^i + z_k^i\big) +  2\min\big\{y_k^i,z_k^i\big\} \big| \nonumber\\
    &\le \big| F\left(y_k - z_k\right)^i - \big(y_k^i + z_k^i\big) \big| +2\big| \min\big\{y_k^i,z_k^i\big\} \big| \nonumber\\
    &\le C r_k^\gamma+\dfrac{4-2{{\tau}}}{{{\tau}}}{r_k}    
    \leq \left(C+\dfrac{4-2{{\tau}}}{{{\tau}}}\right){r_k},
\label{form637}
\end{align}
where we have used Assumption~\ref{assumption1}\,(iii) with the computational error $\e_k^{(1)}=C r_k^\gamma$ and~\eqref{compl-err-i} to conclude.

\medskip

(iii) Thanks to Lemma~\ref{lem:H_global}, we know that $H$ is globally defined and hence 
\[
H(y_k)-F(y_k-H(y_k))+y_k=0 \quad \text{ for every $k \in \mathbb N$.}
\]
It follows
    \begin{align*}
        O( r_k^\gamma) \ge \big\|z_k-F(y_k-z_k)+y_k \big\| &= \big\|z_k-F(y_k-z_k)+y_k-(H(y_k)-F(y_k-H(y_k))+y_k) \big\|\\
        &= \big\|(F+I)(y_k-H(y_k))-(F+I)(y_k-z_k) \big\|.
    \end{align*}
Since $F-I$ is a uniform $P$-map with modulus $h(t)\ge t^\beta$, the mapping $F+I$ has the same property, which gives us
\begin{align*}
    \|z_k-H(y_k)\|^\beta \le & ~ \max\limits_{1\le i\le d}\left((F+I)(y_k-H(y_k))^i-(F+I)(y_k-z_k)^i\right)(z^i_k-H(y_k)^i) \\
    \le & ~ O( r_k^\gamma) \|z_k-H(y_k)\|,
\end{align*}
and so $\|z_k-H(y_k)\|\le O(  r_k^{\frac{\gamma}{\beta-1}}).$

\medskip

(iv) If $y_k$ is bounded, then $z_k$ is bounded by (iii). If $y_k$ is not bounded, then there exists a subsequence
(still denoted $y_k$) such that $\|y_k\|\to +\infty$.
From the first estimate in~\eqref{compl-err-i}, we deduce that $z_k\to 0$. Hence $H(y_k)$ is bounded by (iii).
Using the fact that $F-I$ is uniform-$P$ with modulus $h(t)\ge t^\beta$, we have
\begin{eqnarray}
  \|y_k-H(y_k)\|^\beta &\leq& \max\limits_{1\le i\le d}\left((F-I)(y_k-H(y_k))^i-(F-I)(0)^i\right)(y^i_k-H(y_k)^i)\nonumber\\
  &=&  2 \max\limits_{1\le i\le d} (H(y_k)^i-H(0)^i) y_k^i -  (H(y_k)^i-H(0)^i)^2\nonumber\\
  &\leq & 2 \| H(y_k)-H(0)\| \|y_k\|, \label{ineq481}
\end{eqnarray} 
from which we deduce that $y_k$ is bounded (since $\beta >1$), which is a contradiction.
Finally, both $y_k$ and $z_k$ are bounded.

\medskip

(v) To establish the complementarity errors in terms of $H(y_k)$ instead of $z_k$,
we first note that, since $\psi$ is Lipschitz continuous with constant $L\leq |\psi'(0)|$
and thanks to (iii),
we can replace $\tilde{z}= z_k^i/r_k$ in~\eqref{lem222}
by $\hat{z} = H(y_k)^i/r_k$, i.e.,
\begin{eqnarray}\label{ineq690}
   \psi\left(y_k^i / r_k\right)+ \psi\left(H(y_k)^i / r_k\right)\geq \psi(\e_k^{(2)}/r_k ) -L r_k^{-1}\|z_k-H(y_k)\|
\geq  \psi(\e_k^{(2)}/r_k )+  O(  r_k^{\frac{\gamma -\beta +1}{\beta-1}}).&&
\end{eqnarray}
If $\gamma >\beta -1$, then $O(  r_k^{\frac{\gamma -\beta +1}{\beta-1}})\to 0$ as $k\to \infty$. Therefore, up to
take $k$ large enough, the right hand side of the previous inequality is larger than $\tilde{\tau}>0$
(which may be smaller than $\tau$).
It follows that the first inequality in~\eqref{compl-err-i} holds with $H(y_k)^i$, $\tilde{\tau}$
instead of $z_k$, $\tau$ respectively. Hence the first estimate in~\eqref{complem-iv-1}.

We reproduce the proof of the second inequality in~\eqref{compl-err-i} with
$$
\beta_k:=  \psi(\e_k^{(2)}/r_k )+  O(  r_k^{\frac{\gamma -\beta +1}{\beta-1}}).
$$
Inequality~\eqref{betk-2} now reads
\begin{eqnarray*}
  &&  (1-\beta_k) =\frac{1-\psi(\e_k^{(2)}/r_k)}{\e_k^{(2)}/r_k}C r_k^{\sigma -1}
  + O(  r_k^\frac{\gamma -\beta +1}{\beta-1})
\leq C (-\psi'(0)) r_k^{\sigma -1}+ O( r_k^\frac{\gamma -\beta +1}{\beta-1} ),
\end{eqnarray*}
and~\eqref{betk-1} becomes
\begin{eqnarray*}
  y_k^i H(y_k)^i \le \frac{2 - \tilde{\tau}}{\tilde{\tau}} {r^2_k}
  +\frac{C |\psi'(0)|}{\tilde{\tau}}  \left(y_k^i+z_k^i\right) {r^{\sigma}_k}
  +\frac{1}{\tilde{\tau}}\left(y_k^i+z_k^i\right) O( r_k^\frac{\gamma}{\beta-1} ).
\end{eqnarray*}
Since $y_k^i+z_k^i$ is bounded by (iv), we conclude the second estimate in~\eqref{complem-iv-1}.

\medskip

(vi) Replacing $z_k$ by $H(y_k)$ in~\eqref{form637}
and using $F\left(y_k - H(y_k)\right) - \left(y_k + H(y_k)\right) = 0$
as well as $\min\{y_k^i,H(y_k)^i\}\le O\left(r_k\right)$ by (v), we conclude.\hfill$\Box$


\subsection[]{Proof of \Cref{merged_estimates}}\label{proof_merged_estimates}
Let $z_* = H(y_*)$ and $\widetilde{x}_k\coloneqq y_k-H(y_k)$.
Using that $0\le y_* \perp z_*\ge 0$ and $y_k,z_k\ge 0$,~\Cref{1&2eq-est-bis} yields, for all $1\leq i\leq d$,
\begin{eqnarray}
\left( z^i_* - H(y_k)^i\right) \left( y^i_* - y_k^i\right) &=&  z^i_*y^i_* -z^i_* y_k^i - H(y_k)^iy^i_* + H(y_k)^iy_k^i \nonumber\\
&=&  -z^i_* y_k^i -y^i_* z_k^i + y^i_*(z_k^i -  H(y_k)^i) + y_k^iH(y_k)^i\nonumber\\
&\leq &
\|y_*\| \|z_k-H(y_k)\| + y_k^iH(y_k)^i\nonumber\\
&\leq & 
O( r_k^{\frac{\gamma}{\beta-1}} ) + O( r_k^{\tilde{\sigma}} ) =  O( r_k^{\tilde{\sigma}} ),
\label{upper_max}
\end{eqnarray}
since $\tilde{\sigma}=\min \{ 2, \sigma , \frac{\gamma}{\beta-1}\}$.
Using the uniform $P$-property of $F-I$, it follows
\begin{eqnarray*}
\left\|{x_*} -\widetilde{x}_k\right\|^\beta \le h\left(\left\|{x_*} -\widetilde{x}_k\right\|\right)
&\le &\
\max_{1\le i\le d} \left( (F-I)(x_*)^i - (F-I)\left(\widetilde{x}_k\right)^i\right)\left( {x^i_*} - {\widetilde{x}^i_k}\right)\\
&=&
2\max_{1\le i\le d} \left( z_*^i-H(y_k)^i\right) \left(  y^i_* - y^i_k \right)
   - \left( z_*^i -H(y_k)^i\right)^2\\
&\leq &
 O( r_k^{\tilde{\sigma}} ).   
\end{eqnarray*}
We have
\begin{align*}
        \left\|{x_*} -\widetilde{x}_k\right\|^2 = \left\|y_* - z_* - y_k + H(y_k)\right\|^2
        = \left\|y_* -y_k\right\|^2 + \left\|z_* - H(y_k)\right\|^2 - 2\left\langle y_* -y_k,z_* - H(y_k) \right\rangle. 
\end{align*}
Therefore, taking into account the previous estimate, we obtain
\begin{eqnarray*}
  \left\|y_* -y_k\right\|,  \left\|z_* - H(y_k)\right\| \leq  O( r_k^{\tilde{\sigma}/\beta} ) +  O( r_k^{\tilde{\sigma}/2} )
  =  O( r_k^{\alpha} )  \quad \text{ with } \quad \alpha:= \min \left\{ \frac{\tilde{\sigma}}{\beta} , \frac{\tilde{\sigma}}{2} \right\}.
\end{eqnarray*}
Finally, using~\Cref{1&2eq-est-bis}\,(iii), we get
\begin{align*}
  \|x_*-x_k\| = \|x_*-\widetilde{x}_k\|+\|\widetilde{x}_k-x_k\| = & ~\|x_*-\widetilde{x}_k\|+\|z_k-H(y_k)\| \\
  = & ~ O( r_k^{\tilde{\sigma}/\beta} ) + O( r_k^{\gamma/(\beta-1)} ) = O( r_k^{\tilde{\sigma}/\beta} ).
\end{align*}
Letting $k\to\infty$ and noting that $r_k\to 0$, we get the convergence of $y_k,z_k$ and $x_k$ to $y_*,H(y_*)$ and $x_*$, respectively.
\hfill$\Box$
\subsection[]{Guaranteeing assumptions for reformulating  \ref{NAVE} in specific cases of \Cref{ss3.1}}\label{guarantee_sec3.1}
Suppose that $\mu>\lambda\ge 0$, i.e., $\mu^j>\lambda^j\ge 0$ for all $j=1,\cdots,d$. Then
\begin{align*}
    (\mu - \lambda)^{-1} (\nabla \mathcal{L}(x) + 2\lambda x) = \text{diag}\left(\left(\mu^j-\lambda^j\right)^{-1}\right)(A^TAx-A^Tb)+2\text{diag}\left(\lambda^j\left(\mu^j-\lambda^j\right)^{-1}\right)x.
\end{align*}
Thus
\begin{align*}
    \nabla \left[(\mu - \lambda)^{-1} (\nabla \mathcal{L}(x) + 2\lambda x)\right] = \text{diag}\left(\left(\mu^j-\lambda^j\right)^{-1}\right)A^TA+2\text{diag}\left(\lambda^j\left(\mu^j-\lambda^j\right)^{-1}\right).
\end{align*}
Since $2\text{diag}\left(\lambda^j\left(\mu^j-\lambda^j\right)^{-1}\right)$ and $A^TA$ are positive semi-definite and $\text{diag}\left(\left(\mu^j-\lambda^j\right)^{-1}\right)$ is a
positive diagonal matrix, we thus have the positive semi-definition of the right hand side. This proves the monotonicity of the operator $(\mu - \lambda)^{-1} (\nabla \mathcal{L} + 2\lambda I)$. 
Here, we apply the result established in \cite{NR2022}, which states that the product of a positive diagonal matrix and a symmetric positive semi-definite matrix remains positive semi-definite.
\end{document}